\documentclass[preprint,12pt]{elsarticle}

\usepackage{amssymb}
\usepackage{amsthm}
\usepackage{amsmath}
\usepackage{fancybox}

\newtheorem{thrm}{Theorem}

\newtheorem{thm}{Theorem}
\newtheorem{lem}{Lemma}
\newdefinition{rem}{Remark}

\catcode`\@=11  \@addtoreset{rem}{section}  \catcode`\@=12
\newproof{pf}{{\sc Proof}}
\newproof{pfirst}{Proof of Theorem~\protect\ref{t:first}}
\newproof{pweak}{Proof of Theorem~\protect\ref{t:weak}}
\newproof{prtB}{Proof of Theorem~\protect\ref{t:B}}
\newproof{prLdeux}{Proof of Lemma~\protect\ref{l:deux}}
\newproof{prLtroix}{Proof of Lemma~\protect\ref{l:troix}}
\newproof{prEtwo}{Proof of Lemma~\protect\ref{l:Etwo}}
\newproof{prEone}{Proof of Lemma~\protect\ref{l:Eone}}
\newproof{prtE}{Proof of Theorem~\protect\ref{t:E}}
\newcommand{\wto}{\stackrel{w\,}{\to}}
\newcommand{\eq}[1]{\mbox{\rm(\ref{#1})}}
\newcommand{\ov}[1]{\overline{#1}}
\newcommand{\even}{\mbox{\scriptsize\rm ev\,}}
\newcommand{\odd}{\mbox{\scriptsize\rm od\,}}

\catcode`\@=11  \@addtoreset{equation}{section}  \catcode`\@=12
\newcommand{\LC}{\mathcal{L}}   \newcommand{\RC}{\mathcal{R}}
\newcommand{\PC}{\mathcal{P}}

\journal{arXiv}

\begin{document}

\begin{frontmatter}

\title{Maps of several variables of finite total variation\\
       and Helly-type selection principles\tnoteref{t1}}
\tnotetext[t1]{The work of V.\,V.~Chistyakov was supported by the State University
Higher School of Economics Grant No.~09-08-0012 (Priority thematics
``In\-ves\-ti\-ga\-tion of Function Spaces'').}

\author[hsenn]{Vyacheslav V.~Chistyakov\corref{cor1}}
\ead{czeslaw@mail.ru}

\author[hsenn]{Yuliya V.~Tretyachenko}
\ead{tretyachenko\_y\_v@mail.ru}

\cortext[cor1]{Corresponding author; tel.: +7 831 4169649; fax: +7 831 4169650.}

\address[hsenn]{Department of Applied Mathematics and Informatics,
State University Higher School \endgraf of Economics, Bol'shaya Pech{\"e}rskaya
Street 25/12, Nizhny Novgorod 603155, Russia}

\begin{abstract}
Given two points $a=(a_1,\dots,a_n)$ and $b=(b_1,\dots,b_n)$ from $\mathbb{R}^n$
with $a<b$ componentwise and a map $f$ from the rectangle $I_a^b=[a_1,b_1]\times\cdots\times[a_n,b_n]$
into a metric semigroup $M=(M,d,+)$, we study properties of the {\em total variation\/}
$\mbox{TV}(f,I_a^b)$ of $f$ on $I_a^b$ introduced by the first author in
[V.\,V.~Chistyakov, A selection principle for mappings of bounded variation of several variables,
in: Real Analysis Exchange 27th Summer Symposium, 2003, 217--222], which extends the classical
notion of C.~Jordan's total variation ($n=1$) and the corresponding notions in the sense of
[T.\,H.~Hildebrandt, Introduction to the Theory of Integration, Academic Press, 1963] ($n=2$)
and [A.\,S.~Leonov, On the total variation for functions of several variables and a
multidimensional analog of Helly's selection principle, Math.\ Notes 63 (1998) 61--71]
($n\in\mathbb{N}$) for real valued functions of $n$ variables. The following Helly-type
pointwise selection principle is proved: {\em If a sequence $\{f_j\}_{j\in\mathbb{N}}$ of maps
from $I_a^b$ into $M$ is such that the closure in $M$ of the set $\{f_j(x)\}_{j\in\mathbb{N}}$
is compact for each $x\in I_a^b$ and $C\equiv\sup_{j\in\mathbb{N}}\mbox{\rm TV}(f_j,I_a^b)$ is finite,
then there exists a subsequence of $\{f_j\}_{j\in\mathbb{N}}$, which converges pointwise on
$I_a^b$ to a map $f$ such that\/ \mbox{$\mbox{\rm TV}(f,I_a^b)\le C$}.}
A~variant of this result is established concerning the weak pointwise convergence when
values of maps lie in a reflexive Banach space $(M,\|\cdot\|)$ with separable dual~$M^*$.
\end{abstract}


\end{frontmatter}

\section{Introduction} \label{s:intro}

The classical Helly selection principle (\cite{helly}) states that
{\em a bounded sequence of real valued functions on the closed
interval, which is of uniformly bounded\/ {\rm(}Jordan\/{\rm)} variation,
contains a pointwise convergent subsequence whose limit is a function
of bounded variation}. This theorem and its recent generalizations for real valued
functions and metric space valued maps of one real variable
(\cite{jmaa00,msb98,Sovae,jmaa05,sibam06,jmaa08,waterman})
have numerous applications in different branches of Analysis (e.g.,
\cite{Barbu,Sovae,GNW,hildebrandt,moreau} and references therein).

Extensions of the Helly theorem for functions and maps of several real variables
heavily depend upon notions of (bounded) variation used for these maps, which
generalize different aspects of the classical Jordan variation of univariate
functions and which are known to be quite numerous in the literature (e.g.,
\cite{AFP,Billi,clarkson,Giusti,hildebrandt,ivanov,picone,vitali,vitushkin,ziemer},
and these references are far from being exhaustive on the topic).
Under some approaches to the multidimensional variation (\cite{ambrosio,Billi,nadir})
involving integration procedures Helly-type theorems are rather concerned
with the {\em almost everywhere\/} convergence of extracted subsequences,
and no stronger convergence can be expected in this case, but this convergence
is far too weak for certain applications (such as those from \cite{Sovae}).
On the other hand, there are definitions of the notion of variation for
{\em real valued\/} functions of several variables (\cite{hildebrandt,leonov}),
which go back to Vitali \cite{vitali}, Hardy \cite{hardy} and Krause \cite{adams,clarkson},
such that a complete analogue of the Helly theorem holds with respect to the
{\em pointwise\/} convergence of extracted subsequences. These counterparts of
Helly's theorem are based on the notion of a (totally) {\em monotone\/} real valued function of
several variables \cite{bochner,hildebrandt,young} and an appropriate generalization
of Jordan's decomposition theorem when a function of bounded variation is
represented as the difference of two monotone functions.

The aim of this paper is twofold. First, we study properties of the
{\em total variation\/} of metric semigroup valued maps of several variables
in the approach of Vitali, Hardy and Krause
introduced by the first author in \cite{Opava}, which extends the classical
notion of Jordan's total variation for maps of one variable
and the notions of the total variation in the sense of Hildebrandt \cite{hildebrandt}
for real valued functions of two variables and Leonov \cite{leonov}
for real valued functions of any finite number of variables.
Second, we present two variants of a Helly-type {\em pointwise\/} selection principle
for metric semigroup valued maps and maps with values in a reflexive separable Banach space.
The main difficulty that we overcome is that for metric semigroup valued maps there is
no counterpart of Jordan's decomposition theorem, and we have to develop a
completely different technique, whose two-dimesional variant is given in \cite{austral}.

\smallbreak
The paper is organized as follows. In Section~\ref{s:defs} we present necessary
definitions and our two main results, Theorems~\ref{t:first} and~\ref{t:weak}.
In order to get to the proofs of these results as quick as possible, in Section~\ref{s:ingred}
we collect all main ingredients and auxiliary known facts needed for their proofs.
In Section~\ref{s:proofs} we prove Theorems~\ref{t:first} and~\ref{t:weak}.
The remaining Sections \ref{s:prB}--\ref{s:prE} contain proofs of the results
exposed in Section~\ref{s:ingred} and used in the proofs of the main~theorems.

\section{Definitions and main results} \label{s:defs}

Throughout the paper we adopt and follow the Vitali-Hardy-Krause approach
to the notion of variation for maps of several variables in the multiindex
notation initiated in \cite{lexington,Opava} and developed in detail in~\cite{na05}
(equivalent approaches in different notation for real functions can be
found in \cite{lenze,leonov}).

Let $\mathbb{N}$ and $\mathbb{N}_0$ stand for the sets of positive and nonnegative
integers, respectively, and $n\in\mathbb{N}$. Given $x,y\in\mathbb{R}^n$, we write
$x=(x_1,\dots,x_n)=(x_i:i\in\{1,\dots,n\})$ for the coordinate representation of $x$,
and set $x+y=(x_1+y_1,\dots,x_n+y_n)$, and $x-y$ is defined similarly.
The inequality $x<y$ will be understood componentwise, i.e.,  $x_i<y_i$ for all
$i\in\{1,\dots,n\}$, and a similar meaning applies to $x=y$, $x\le y$, $y\ge x$
and $y>x$. If $x<y$ or $x\le y$, we denote by $I_x^y$ the rectangle
$\prod_{i=1}^n[x_i,y_i]=[x_1,y_1]\times\cdots\times[x_n,y_n]$.
Elements of the set $\mathbb{N}_0^n$ are as usual said to be {\em multiindices\/}
and denoted by Greek letters and, given $\theta=(\theta_1,\dots,\theta_n)\in\mathbb{N}_0^n$
and $x\in\mathbb{R}^n$, we set $|\theta|=\theta_1+\cdots+\theta_n$ (the order
of $\theta$) and $\theta x=(\theta_1x_1,\dots,\theta_nx_n)$.
The $n$-dimensional multiindices $0_n=(0,\dots,0)$ and $1_n=(1,\dots,1)$ will be denoted
simply by $0$ and $1$, respectively (actually, the dimension of $0$ and $1$ will be clear
from the context). We also put
$\mathcal{E}(n)=\{\theta\in\mathbb{N}_0^n:\mbox{$\theta\le1$ and $|\theta|$ is even}\}$
(the set of `even' multiindices) and
$\mathcal{O}(n)=\{\theta\in\mathbb{N}_0^n:\mbox{$\theta\le1$ and $|\theta|$ is odd}\}$
(the set of `odd' multiindices).
For elements from the set $\mathcal{A}(n)=\{\alpha\in\mathbb{N}_0^n:0\ne\alpha\le1\}$
we simply write $0\ne\alpha\le1$.

The domain of (almost) all maps under consideration will be a rectangle $I_a^b$
with fixed $a,b\in\mathbb{R}^n$, $a<b$, called the {\em basic rectangle}.
The range of maps will be a {\em metric semigroup\/} $(M,d,+)$, i.e., $(M,d)$
is a metric space, $(M,+)$ is an Abelian semigroup with the operation of addition~$+$,
and $d$ is translation invariant: $d(u,v)=d(u+w,v+w)$ for all $u,v,w\in M$.
A nontrivial example of a metric semigroup is as follows (\cite{deblasi,radstrom}):
Let $(X,\|\cdot\|)$ be a real normed space and $M$ be the family of all nonempty
closed bounded convex subsets of $X$ equipped with the Hausdorff metric $d$ given by
$d(U,V)=\max\{\mbox{e}(U,V),\mbox{e}(V,U)\}$, where $U,V\in M$ and
$\mbox{e}(U,V)=\sup_{u\in U}\inf_{v\in V}\|u-v\|$. Given $U,V\in M$,
defining $U\oplus V$ as the closure in $X$ of the Minkowski sum
$\{u+v:u\in U,v\in V\}$ we find that the triple $(M,d,\oplus)$ is a metric semigroup.

Given $f:I_a^b\to(M,d,+)$, we define the {\em Vitali-type $n$-th mixed `difference'}
of $f$ on a subrectangle $I_x^y\subset I_a^b$, where $x,y\in I_a^b$ and $x<y$,
by (cf.~\cite{Opava})
  \begin{equation} \label{e:mdn}
\mbox{md}_n(f,I_x^y)=d\Bigl(\sum_{\theta\in\mathcal{E}(n)}f\bigl(x+\theta\,(y-x)\bigr),
\sum_{\eta\in\mathcal{O}(n)}f\bigl(x+\eta\,(y-x)\bigr)\Bigr).
  \end{equation}

For example, for the first three dimensions we have: if $n=1$, then $\mathcal{E}(1)=\{0\}$
and $\mathcal{O}(1)=\{1\}$, and so, $\mbox{md}_1(f,I_x^y)=d(f(x),f(y))$;
if $n=2$, then $\mathcal{E}(2)=\{(0,0),(1,1)\}$ and $\mathcal{O}(2)=\{(0,1),(1,0)\}$,
and so,
  $$\mbox{md}_2(f,I_x^y)=d\bigl(f(x_1,x_2)+f(y_1,y_2),f(x_1,y_2)+f(y_1,x_2)\bigr);$$
if $n=3$, then $\mathcal{E}(3)\!=\!\{(0,0,0),(1,1,0),(1,0,1),(0,1,1)\}$ and
$\mathcal{O}(3)\!=\!\{(1,1,1),$ $(1,0,0),(0,1,0),(0,0,1)\}$, and so,
  \begin{eqnarray*}
\mbox{md}_3(f,I_x^y)\!&=\!&d\Bigl(f(x_1,x_2,x_3)+f(y_1,y_2,x_3)+f(y_1,x_2,y_3)
  +f(x_1,y_2,y_3),\\[-2pt]
\!&\!&\quad\,\,f(y_1,y_2,y_3)+f(y_1,x_2,x_3)+f(x_1,y_2,x_3)+f(x_1,x_2,y_3)\Bigr)
  \end{eqnarray*}
(one may draw corresponding pictures to see the points where $f$ is evaluated
at the left and right hand places of $d(\mbox{`left'},\mbox{`right'})$\,).

\begin{rem} \label{remA}
Formally, the value $\mbox{md}_n(f,I_x^y)$ from \eq{e:mdn} is defined for $x<y$.
Now if $x,y\in I_a^b$, $x\le y$ and $x\not<y$, then the
right-hand side in \eq{e:mdn} is equal to zero for any map $f\!:\!I_a^b\!\to\! M$.
In fact, if $x_i\!=\!y_i$ for some $i\!\in\!\{1,\dots,n\}$, then
  $$
\sum_{\theta\in\mathcal{E}(n)}f(x+\theta(y-x))=\sum_{\ov\theta\in\mathcal{O}(n)}f(x+\ov\theta(y-x)).
  $$
In order to see this, given $\theta=(\theta_1,\dots,\theta_n)\in\mathcal{E}(n)$, we set
  $\ov\theta=(\ov\theta_1,\dots,\ov\theta_n)=(\theta_1,\dots,\theta_{i-1},1-\theta_i,
    \theta_{i+1},\dots,\theta_n)$
and note that $\ov\theta\in\mathcal{O}(n)$ and, moreover, the map $\theta\mapsto\ov\theta$
is a bijection between $\mathcal{E}(n)$ and $\mathcal{O}(n)$. It remains to take into
account that $x+\theta(y-x)=x+\ov\theta(y-x)$ for all $\theta\in\mathcal{E}(n)$, because
  $$x_i+\theta_i(y_i-x_i)=x_i=x_i+(1-\theta_i)(y_i-x_i)=x_i+\ov\theta_i(y_i-x_i).$$
\end{rem}

The {\em Vitali-type $n$-th variation\/} (\cite{na05,leonov,vitali})
of $f:I_a^b\to M$ is defined by
  \begin{equation} \label{e:Vn}
V_n(f,I_a^b)=\sup_{\cal P}\sum_{1\le\sigma\le\kappa}
\mbox{md}_n\bigl(f,I_{x[\sigma-1]}^{x[\sigma]}\bigr),
  \end{equation}
the supremum being taken over all multiindices $\kappa\in\mathbb{N}^n$ and
all {\em net partitions\/} of $I_a^b$ of the form
$\mathcal{P}=\{x[\sigma]\}_{\sigma=0}^\kappa$, where
points $x[\sigma]=(x_1(\sigma_1),\dots,x_n(\sigma_n))$ from $I_a^b$ are
indexed by $\sigma=(\sigma_1,\dots,\sigma_n)\in\mathbb{N}_0^n$ with
$\sigma\le\kappa$ and satisfy the conditions: $x[0]=a$, $x[\kappa]=b$ and
$x[\sigma-1]<x[\sigma]$ for all $1\le\sigma\le\kappa$ (in other words, a net partition
$\mathcal P$ is the Cartesian product of ordinary partitions of closed intervals
$[a_i,b_i]$, $i=1,\dots,n$). Note that all rectangles $I_{x[\sigma-1]}^{x[\sigma]}$
of a net partition are non-degenerated, non-overlapping and their union is $I_a^b$.

In order to define the notion of the total variation of a map $f:I_a^b\to M$ we need
the notion of variation of $f$ of order less than~$n$. Following \cite{na05},
we define the {\em truncation of a point $x\in\mathbb{R}^n$ by a multiindex
$0\ne\alpha\le1$} by $x\lfloor\alpha=(x_i\,:\,i\in\{1,\dots,n\},\,\alpha_i=1)$,
and set $I_a^b\lfloor\alpha=I_{a\lfloor\alpha}^{b\lfloor\alpha}$.
Clearly, $x\lfloor1=x$ and $I_a^b\lfloor1=I_a^b$, and if $x\in I_a^b$, then
$x\lfloor\alpha\in I_a^b\lfloor\alpha$.
For example, if $x=(x_1,x_2,x_3,x_4)$ and $\alpha=(1,0,0,1)$,
we have $x\lfloor\alpha=(x_1,x_4)$ and
$I_a^b\lfloor\alpha=[a_1,b_1]\times[a_4,b_4]$.
Given $f:I_a^b\to M$ and $z\in I_a^b$, we define the {\em truncated map\/}
$f_\alpha^z:I_a^b\lfloor\alpha\to M$ {\em with the base at $z$} by
$f_\alpha^z(x\lfloor\alpha)=f\bigl(z+\alpha\,(x-z)\bigr)$ for all $x\in I_a^b$.
It follows that $f_\alpha^z$ depends only on $|\alpha|$ variables $x_i\in[a_i,b_i]$,
for which $\alpha_i=1$, and the other variables remain fixed and equal to $z_j$
when $\alpha_j=0$. In the above example we get
$f_\alpha^z(x_1,x_4)=f_\alpha^z(x\lfloor\alpha)=f(x_1,z_2,z_3,x_4)$ for
$(x_1,x_4)\in[a_1,b_1]\times[a_4,b_4]$.

Now, given $f:I_a^b\to M$ and $0\ne\alpha\le1$, the function
$f_\alpha^a:I_a^b\lfloor\alpha\to M$ with the base at $z=a$
depends only on $|\alpha|$ variables, and so, making use of the definitions
(\ref{e:Vn}) and (\ref{e:mdn}) with $n$ replaced by $|\alpha|$, $f$ replaced
by $f_\alpha^a$ and $I_a^b$ replaced by $I_a^b\lfloor\alpha$, we get the notion
of the ({\em Hardy-Krause-type\/} \cite{adams,clarkson,hardy})
{\em$|\alpha|$-th var\-i\-a\-tion\/} of $f$, which is denoted by
$V_{|\alpha|}(f_\alpha^a,I_a^b\lfloor\alpha)$.

The {\em total variation\/} of $f:I_a^b\to M$ in the sense of Hildebrandt
(\cite{dan03,smz05}, \cite[III.6.3]{hildebrandt}, \cite{idczak} if $n=2$)
and Leonov (\cite{lexington,Opava,na05,leonov} if $n\in\mathbb{N}$) is defined by
  \begin{equation} \label{e:TV}
\mbox{TV}(f,I_a^b)=\sum_{0\ne\alpha\le1}V_{|\alpha|}(f_\alpha^a,I_a^b\lfloor\alpha)
\equiv\sum_{i=1}^n\sum_{\alpha\le1,\,|\alpha|=i}V_i(f_\alpha^a,I_a^b\lfloor\alpha),
  \end{equation}
the summations here and throughout the paper being taken over
{\em $n$-di\-men\-sion\-al\/} multiindices in the ranges specified
under the summation sign.

For the first three dimensions $n=1,2,3$ we have, respectively,
  \begin{eqnarray*}
\mbox{TV}(f,I_a^b)&=&V_a^b(f),\quad\mbox{the usual Jordan variation on the interval $[a,b]$},\\
\mbox{TV}(f,I_a^b)&=&V_{a_1}^{b_1}\bigl(f(\,\cdot\,,a_2)\bigr)
  +V_{a_2}^{b_2}\bigl(f(a_1,\,\cdot\,)\bigr)+V_2(f,I_{a_1,a_2}^{b_1,b_2}),\\
\mbox{TV}(f,I_a^b)&=&V_{a_1}^{b_1}\bigl(f(\,\cdot\,,a_2,a_3)\bigr)
  +V_{a_2}^{b_2}\bigl(f(a_1,\,\cdot\,,a_3)\bigr)
  +V_{a_3}^{b_3}\bigl(f(a_1,a_2,\,\cdot\,)\bigr)\\[1pt]
&&\quad+V_2\bigl(f(\,\cdot\,,\,\cdot\,,a_3),I_{a_1,a_2}^{b_1,b_2}\bigr)
  +V_2\bigl(f(\,\cdot\,,a_2,\,\cdot\,),I_{a_1,a_3}^{b_1,b_3}\bigr)\\[1pt]
&&\quad+V_2\bigl(f(a_1,\,\cdot\,,\,\cdot\,),I_{a_2,a_3}^{b_2,b_3}\bigr)
  +V_3(f,I_{a_1,a_2,a_3}^{b_1,b_2,b_3}).
  \end{eqnarray*}

We denote by $\mbox{BV}(I_a^b;M)$ the space of all maps $f:I_a^b\to M$ of {\em finite\/}
(or bounded) {\em total variation\/} (\ref{e:TV}).

Recall that a sequence $\{f_j\}\equiv\{f_j\}_{j\in\mathbb{N}}$ of maps from $I_a^b$
into $M$ is said: (a) to {\em converge pointwise\/} on $I_a^b$ to a map $f:I_a^b\to M$
if $d(f_j(x),f(x))\to0$ as $j\to\infty$ for all $x\in I_a^b$;
(b) to be {\em pointwise precompact\/} (on $I_a^b$) provided the closure in $M$
of the set $\{f_j(x)\}_{j\in\mathbb{N}}$ is compact for all $x\in I_a^b$.

Our first main result, to be proved in Section~\ref{s:proofs}, is the following Helly-type
pointwise selection principle in the space $\mbox{BV}(I_a^b;M)$:

\begin{thm} \label{t:first}
A pointwise precompact sequence $\{f_j\}$ of maps from the rectangle $I_a^b$ into
a metric semigroup $(M,d,+)$ such that
  \begin{equation} \label{e:ubv}
C\equiv\sup_{j\in\mathbb{N}}\mbox{\rm TV}(f_j,I_a^b)\quad\mbox{\,is \,\,finite}
  \end{equation}
contains a subsequence which converges pointwise on $I_a^b$ to a map
$f\!\in\!\mbox{\rm BV}(I_a^b;M)$ such that $\mbox{\rm TV}(f,I_a^b)\le C$.
\end{thm}

This result was announced in \cite{Opava}. It contains as particular cases
the results of \cite[III.6.5]{hildebrandt} and \cite{idczak} ($n=2$ and $M=\mathbb{R}$),
\cite{leonov} ($n\in\mathbb{N}$ and $M=\mathbb{R}$)
and \cite{austral} ($n=2$ and $M$ is a metric semigroup).

\smallbreak
Our second main result (Theorem~\ref{t:weak} below) is concerned with a weak analogue
of Theorem~\ref{t:first} taking into account certain specific features when the values
of maps under consideration lie in a reflexive separable Banach space.

Let $(M,\|\cdot\|)$ be a normed linear space over the field $\mathbb{K}=\mathbb{R}$
or $\mathbb{C}$ and $M^*$ be its {\em dual}, i.e., $M^*=\mbox{\rm L}(M;\mathbb{K})$,
the space of all continuous linear functionals on $M$. It is well-known that $M^*$
is a Banach space under the norm $\|u^*\|=\sup\{|u^*(u)|:\mbox{$u\in M$ and $\|u\|\le1$}\}$,
$u^*\in M^*$. The natural duality between $M$ and $M^*$ is determined by the bilinear
functional $\langle\cdot,\cdot\rangle:M\times M^*\to\mathbb{K}$ defined by
$\langle u,u^*\rangle=u^*(u)$ for all $u\in M$ and $u^*\in M^*$, so that
$|\langle u,u^*\rangle|\le\|u\|\cdot\|u^*\|$, where $|\cdot|$ is the absolute value
in $\mathbb{K}$. Recall that a sequence $\{u_j\}\subset M$ {\em converges weakly in $M$}
to an element $u\in M$ (in symbols, $u_j\wto u$ in $M$) if
$\langle u_j,u^*\rangle\to\langle u,u^*\rangle$ in $\mathbb{K}$ as $j\to\infty$
for all $u^*\in M^*$; if this is the case then it is known that
$\|u\|\le\liminf_{j\to\infty}\|u_j\|$.

Since a normed linear space $(M,\|\cdot\|)$ is a metric semigroup, the notions
of the Vitali-type $n$-th variation, $|\alpha|$-th variation for $0\ne\alpha\le1$
and the total variation of a map $f:I_a^b\to M$ are introduced as above with
respect to the induced metric $d(u,v)=\|u-v\|$, $u,v\in M$.

\begin{thm} \label{t:weak}
Suppose $(M,\|\cdot\|)$ is a reflexive separable Banach space with separable
dual $M^*$ and $\{f_j\}$ is a sequence of maps from $I_a^b$ into~$M$.
If $\{f_j\}$ satisfies condition\/ \eq{e:ubv} from Theorem\/~{\rm\ref{t:first}} and
  \begin{equation} \label{e:pou}
c(x)\equiv\sup_{j\in\mathbb{N}}\|f_j(x)\|\quad\mbox{\,is \,finite \,for \,all $x\in I_a^b$,}
  \end{equation}
then there exists a subsequence of $\{f_j\}$, again denoted by $\{f_j\}$, and a map
$f\in\mbox{\rm BV}(I_a^b;M)$ satisfying $\mbox{\rm TV}(f,I_a^b)\le C$ such that
  \begin{equation} \label{e:wec}
f_j(x)\wto f(x)\quad\mbox{in \,$M$ \,for \,all \,$x\in I_a^b$.}
  \end{equation}
\end{thm}

This theorem will be proved in Section~\ref{s:proofs}. It is an extension of
a weak selection principle from \cite[Chapter~1, Theorem~3.5]{Barbu} for maps
of bounded Jordan variation of one real variable. More comments and remarks on
Theorems \ref{t:first} and \ref{t:weak} can be found in Section~\ref{s:proofs}.

\section{Properties of mixed differences and the total variation} \label{s:ingred}

In this section we collect main ingredients of the proof of Theorem~\ref{t:first}.
These are relations between mixed differences of all orders and properties
of the total variation (\ref{e:TV}). For real valued functions of $n$ variables the
main properties of mixed differences of all orders were elaborated in
\cite{adams,monatsh,na05,clarkson,hildebrandt,leonov,vitali} and for metric semigroup
valued maps of two variables---in \cite{austral,dan03,smz05,picone}.
For our purposes we need their variants in the multiindex notation, as presented
in \cite{na05} with $M=\mathbb{R}$, for maps of $n$ variables with values in a
metric semigroup.

First, we recall several definitions and results for real valued functions.
A function $g:I_a^b\to\mathbb{R}$ is said to be {\em totally monotone\/}
(cf., e.g., \cite[Part~II, Section~3]{na05}, \cite{leonov}) if, given
$0\ne\alpha\le1$ and $x,y\in I_a^b$ with $x\le y$, we have:
  \begin{equation} \label{e:mon}
(-1)^{|\alpha|}\sum_{0\le\theta\le\alpha}(-1)^{|\theta|}g\bigl(x+\theta(y-x)\bigr)\ge0.
  \end{equation}
For real valued functions the sum in \eq{e:mon} (with no factor $(-1)^{|\alpha|}$)
is called the {\em $|\alpha|$-th mixed difference\/} (in the sense of Vitali, Hardy
and Krause) of $g_\alpha^x$ on the rectangle $I_x^y\lfloor\alpha$ and denoted by
$\mbox{\rm md}_{|\alpha|}(g_\alpha^x,I_x^y\lfloor\alpha)$ (however, note the difference
with \eq{e:fax} in the general case). In this case the {\em Vitali $n$-th variation\/}
$V_n(g,I_a^b)$ of $g$ on $I_a^b$ is defined as in \eq{e:Vn} with the mixed
difference at the right-hand side of \eq{e:Vn} replaced by
$\bigl|\mbox{\rm md}_n(g,I_{x[\sigma-1]}^{x[\sigma]})\bigr|$. The other definitions
related to the bounded variation context remain the same as above, and so, we keep
the same notation for real valued functions as well.

Denote by $\mbox{Mon}(I_a^b;\mathbb{R})$ the set of all totally monotone real valued
functions on $I_a^b$. It is known (e.g., from the references above) that if
$g\in\mbox{Mon}(I_a^b;\mathbb{R})$, then $g\in\mbox{BV}(I_a^b;\mathbb{R})$,
the value at the left-hand side of \eq{e:mon} is equal to
$V_{|\alpha|}(g_\alpha^x,I_x^y\lfloor\alpha)$, $g(x)\le g(y)$ and
$\mbox{TV}(g,I_x^y)=g(y)-g(x)$ for all $x,y\in I_a^b$ with $x\le y$.

The following Helly selection principle in the class $\mbox{Mon}(I_a^b;\mathbb{R})$
is due to Leonov \cite[Lemma~3]{leonov} (for totally monotone functions of two
variables it was established in \cite[III.6.5]{hildebrandt} and \cite[Theorem~3.1]{idczak}):

\begin{thrm} \label{t:A}
$\!$An infinite uniformly bounded family of totally monotone functions on $I_a^b$
contains a sequence, which converges pointwise on $I_a^b$ to a function from
$\mbox{\rm Mon}(I_a^b;\mathbb{R})$.
\end{thrm}

It was shown in \cite[Corollary~2]{leonov} that the linear space $\mbox{BV}(I_a^b;\mathbb{R})$
equipped with the norm $\|g\|=|g(a)|+\mbox{TV}(g,I_a^b)$, $g\in\mbox{BV}(I_a^b;\mathbb{R})$,
is a Banach space. This assertion was refined in \cite[Part~I, Theorem~1]{na05}:
the space $\mbox{BV}(I_a^b;\mathbb{R})$ is a Banach algebra with respect to the norm
$\|\cdot\|$, and $\|g\cdot g'\|\le 2^n\|g\|\cdot\|g'\|$ for all $g,g'\in\mbox{BV}(I_a^b;\mathbb{R})$.

Theorem~\ref{t:A} implies Helly's selection principle in the space
$\mbox{BV}(I_a^b;\mathbb{R})$ \cite[Theorem~4]{leonov}: {\em an infinite family
of functions from\/ $\mbox{\rm BV}(I_a^b;\mathbb{R})$, which is bounded under the norm
$\|\cdot\|$, contains a pointwise convergent sequence, whose pointwise limit
belongs to $\mbox{\rm BV}(I_a^b;\mathbb{R})$}. The crucial observation in the proof of
this result is that, given $g\in\mbox{BV}(I_a^b;\mathbb{R})$, if we set
$\nu_g(x)=\mbox{TV}(g,I_a^x)$ and $\pi_g(x)=\nu_g(x)-g(x)$, $x\in I_a^b$,
then (\cite[Theorem~3]{leonov}) $\nu_g,\pi_g\in\mbox{\rm Mon}(I_a^b;\mathbb{R})$,
and Jordan's decomposition holds: $g=\nu_g-\pi_g$ on $I_a^b$; then Theorem~\ref{t:A}
applies to the uniformly bounded families of functions $\{\nu_g\}$ and $\{\pi_g\}$
in the standard way.

\smallbreak
Now let us consider the case of maps $f:I_a^b\to M$ of finite total variation valued
in a metric semigroup $(M,d,+)$. Clearly, there is no counterpart of Jordan's
decomposition for these maps, and so, in order to prove Theorem~\ref{t:first},
we ought to argue in a completely different way. It will be seen later that,
along with Theorem~\ref{t:A}, the following four Theorems \ref{t:B} through
\ref{t:E} are the main ingredients in the proof of Theorem~\ref{t:first}
(in a certain sense replacing the arguments involving Jordan's decomposition).

\begin{thrm} \label{t:B}
If\/ $f\in\mbox{\rm BV}(I_a^b;M)$, $x,y\in I_a^b$ and $x\le y$, then
  $$d(f(x),f(y))\le\sum_{0\ne\alpha\le1}\mbox{\rm md}_{|\alpha|}(f_\alpha^x,I_x^y\lfloor\alpha)
    \le\mbox{\rm TV}(f,I_x^y).$$
\end{thrm}

This theorem will be proved in Section~\ref{s:prB}. It is a generalization of the
well-known property of maps of bounded Jordan variation of one variable and a
counterpart of Leonov's (in)equalities established in \cite[Theorem~2 and Corollary~5]{leonov}
for real valued functions of $n$ variables (cf.\ also \cite[Part~I, Lemma~6 and (3.5)]{na05}).
The inequalities in Theorem~\ref{t:B} are also known for metric semigroup valued
maps of two variables \cite{austral,smz05}. However, in the general case
Theorem~\ref{t:B} needs a different proof as compared to the cases of maps of
one or two variable(s) or $M=\mathbb{R}$.

\begin{thrm} \label{t:C}
If $f\in\mbox{\rm BV}(I_a^b;M)$, $x,y\in I_a^b$, $x\le y$, and $0\ne\gamma\le1$, then
  \begin{eqnarray}
\sum_{0\ne\alpha\le\gamma}V_{|\alpha|}(f_\alpha^x,I_x^y\lfloor\alpha)&=&
  \mbox{\rm TV}\bigl(f,I_x^{x+\gamma(y-x)}\bigr)\nonumber\\
&\le&\mbox{\rm TV}(f,I_a^{x+\gamma(y-x)})-\mbox{\rm TV}(f,I_a^x).
  \label{e:ustar}
  \end{eqnarray}
\end{thrm}

\begin{thrm} \label{t:D}
If $f\in\mbox{\rm BV}(I_a^b;M)$ and if we set $\nu_f(x)=\mbox{\rm TV}(f,I_a^x)$,
$x\in I_a^b$, then for $\nu_f:I_a^b\to\mathbb{R}$, called the {\sl total variation
function of $f$}, we have\/{\rm:} $\nu_f\in\mbox{\rm Mon}(I_a^b;\mathbb{R})$ and\/
$\mbox{\rm TV}(\nu_f,I_a^b)=\mbox{\rm TV}(f,I_a^b)$.
\end{thrm}

These two theorems are extensions of two more properties of the Jordan variation
for maps of one variable; in this case \eq{e:ustar} is actually the equality
known as the {\em additivity\/} of Jordan's variation (e.g., \cite[Theorem~$83_2$]{schwartz}).
On the other hand, Theorem~\ref{t:C} is a counterpart of Chistyakov's inequality
\cite[Part~II, Lemma~8]{na05} and Theorem~\ref{t:D} is a generalization of
Theorem~3 from \cite{leonov} and Corollary~11 from \cite[Part~II]{na05} given
for $M=\mathbb{R}$. For metric semigroup valued maps of two variables cf.\
\cite[inequalities (11), (13) and Theorem~1]{austral}.

The proof of Theorem~\ref{t:C} is identical with the proof of Lemma~8 from
\cite[Part~II]{na05} and the proof of Theorem~\ref{t:D} is identical with
the proofs of Lemma~9 and Corollaries 10 and 11 from \cite[Part~II]{na05}
when $M=\mathbb{R}$, and so, they are omitted.
However, it is to be noted that these proofs rely on
(1)~equality (3.2) from \cite[Part~I, Lemma~5]{na05}, (2) Lemma~7 from \cite[Part~I]{na05},
and (3) the well-known property of the {\em additivity\/} of $|\alpha|$-th variation
$V_{|\alpha|}$ for each \mbox{$0\ne\alpha\le1$} for real valued functions of $n$ variables.
For metric semigroup valued maps assertions (1), (2) and (3) need a proper
interpretation and different, more subtle and hard proofs. Their respective
counterparts are presented below as Lemmas~\ref{l:un}, \ref{l:deux} and \ref{l:troix}.

In the first lemma and throughout the paper we use the following short notations:
given $0\ne\alpha\le1$, the sum over `$\mbox{\rm ev}\,\theta\le\alpha$' denotes
the sum over `$\theta\in\mathcal{E}(n)$ s.t.\ $\theta\le\alpha$', where `s.t.' is
the usual abbreviation for `such that', and a similar convention applies to
the sum over `$\mbox{\rm od}\,\theta\le\alpha$'.

\begin{lem} \label{l:un}
If $f:I_a^b\to M$, $x,y\in I_a^b$, $x\le y$, $z\in I_a^b$ and $0\ne\alpha\le1$, then
  \begin{eqnarray}
\mbox{\rm md}_{|\alpha|}(f_\alpha^z,I_x^y\lfloor\alpha)&=&d\biggl(
  \sum_{\even\theta\le\alpha}\!f\bigl(z+\alpha(x-z)+\theta(y-x)\bigr),\nonumber\\
&&\quad\,\,\,\sum_{\odd\theta\le\alpha}\!f\bigl(z+\alpha(x-z)+\theta(y-x)\bigr)\biggr).\label{e:2.4*}
  \end{eqnarray}
In particular, if $z=a$ or $z=x$, we have, respectively,
  $$\mbox{\rm md}_{|\alpha|}(f_\alpha^a,I_x^y\lfloor\alpha)
    =\mbox{\rm md}_{|\alpha|}(f_\alpha^{a+\alpha(x-a)},I_{a+\alpha(x-a)}^y\lfloor\alpha),$$
  \begin{equation} \label{e:fax}
\mbox{\rm md}_{|\alpha|}(f_\alpha^x,I_x^y\lfloor\alpha)=d\biggl(\,
\sum_{\even\theta\le\alpha}\!f\bigl(x+\theta(y-x)\bigr),
\sum_{\odd\theta\le\alpha}\!f\bigl(x+\theta(y-x)\bigr)\biggr).
  \end{equation}
\end{lem}

The proof of Lemma~\ref{l:un} is the same as in \cite[Part~I, Lemma~5]{na05}
(details are omitted):
we have to note only that $\theta'\in\mathbb{N}_0^{|\alpha|}$ and $|\theta'|$ is even
(odd) if and only if there exists a unique $\theta\in\mathbb{N}_0^n$ s.t.\
$\theta\le\alpha$, $|\theta|$ is even (odd, respectively) and $\theta'=\theta\lfloor\alpha$,
and apply definition (\ref{e:mdn}) where $n$ is replaced by $|\alpha|$.

\smallbreak
Since the total variation \eq{e:TV} is defined via truncated maps with the base
at the point $a$, in our next lemma we present a counterpart of Chistyakov's
equality \cite[Part~I, Lemma~7]{na05} exhibiting the relation between the mixed
difference $\mbox{md}_{|\alpha|}(f_\alpha^x,I_x^y\lfloor\alpha)$ and certain
mixed differences of maps $f_\beta^a$ with the base at $a$ for some $0\ne\beta\le1$.

\begin{lem} \label{l:deux}
If $f:I_a^b\to M$, $0\ne\alpha\le1$ and $x,y\in I_a^b$ with $x\le y$, then
  $$\mbox{\rm md}_{|\alpha|}(f_\alpha^x,I_x^y\lfloor\alpha)\le\sum_{\alpha\le\beta\le1}
    \mbox{\rm md}_{|\beta|}\bigl(f_\beta^a,I_{a+\alpha(x-a)}^{x+\alpha(y-x)}\lfloor\beta\bigr).$$
\end{lem}

The proof of Lemma~\ref{l:deux} is postponed until Section~\ref{s:prdeux}.

\smallbreak
The {\em additivity\/} property of $|\alpha|$-th variation $V_{|\alpha|}$ for each
$0\ne\alpha\le1$, to be proved in Section~\ref{s:prtroix}, is expressed in the following

\begin{lem} \label{l:troix}
Given $f:I_a^b\to M$, $x,y\in I_a^b$ with $x<y$, $z\in I_a^b$ and $0\ne\alpha\le1$,
if $\{x[\sigma]\}_{\sigma=0}^\kappa$ is a net partition of $I_x^y$, then
  \begin{equation} \label{e:addi}
V_{|\alpha|}(f_\alpha^z,I_x^y\lfloor\alpha)=
\sum_{1\lfloor\alpha\le\sigma\lfloor\alpha\le\kappa\lfloor\alpha}
V_{|\alpha|}(f_\alpha^z,I_{x[\sigma-1]}^{x[\sigma]}\lfloor\alpha),
  \end{equation}
where the summation is taken only over those $\sigma_i$ in the range $1\le\sigma_i\le\kappa_i$
with $i\in\{1,\dots,n\}$, for which $\alpha_i=1$.
\end{lem}

The final ingredient in the proof of Theorem~\ref{t:first} is the {\em sequential
lower semicontinuity\/} of the total variation $\mbox{TV}(\cdot,I_a^b)$ to be
established in Section~\ref{s:prE}:

\begin{thrm} \label{t:E}
If a sequence of maps $\{f_j\}$ from $I_a^b$ into $M$ converges pointwise on $I_a^b$
to a map $f:I_a^b\to M$, then
$\mbox{\rm TV}(f,I_a^b)\le\liminf_{j\to\infty}\mbox{\rm TV}(f_j,I_a^b)$.
\end{thrm}

Now we are in a position to prove Theorems \ref{t:first} and \ref{t:weak}.

\section{Proofs of Theorems \ref{t:first} and \ref{t:weak}} \label{s:proofs}

\begin{pfirst}
We divide the proof into four steps for clarity.

1. We apply the induction argument on the dimension $n$ of the basic rectangle
$I_a^b\subset\mathbb{R}^n$. For $n=1$ Theorem~\ref{t:first} was established in
\cite[Theorem~5.1]{msb98} (and refined in \cite[Theorem~1]{jmaa00} and \cite[Theorem~1.3]{Sovae})
in the case when $(M,d)$ is an arbitrary metric space, and for $n=2$ it was proved
in \cite[Theorem~2]{austral}. Now, suppose that $n\ge3$ and Theorem~\ref{t:first}
is already established for domain rectangles of dimension $\le n-1$.

Given $j\in\mathbb{N}$, we let $\nu_j$ be the total variation function of $f_j$
on $I_a^b$, i.e., $\nu_j(x)=\mbox{TV}(f_j,I_a^x)$ for all $x\in I_a^b$.
By Theorem~\ref{t:D} and condition \eq{e:ubv}, the sequence
$\{\nu_j\}\subset\mbox{Mon}(I_a^b;\mathbb{R})$ is uniformly bounded (by $C$),
and so, by Theorem~\ref{t:A}, there exist a subsequence of $\{\nu_j\}$ and the
corresponding subsequence of $\{f_j\}$, again denoted as the whole sequences $\{\nu_j\}$
and $\{f_j\}$, respectively, and a function $\nu\in\mbox{Mon}(I_a^b;\mathbb{R})$ s.t.\
  \begin{equation} \label{e:nujnu}
\lim_{j\to\infty}\nu_j(x)=\nu(x)\quad\mbox{for \,\,all}\quad x\in I_a^b.
  \end{equation}
It is known (\cite{antosik}, \cite[III.5.4]{hildebrandt}, \cite{young}) that the
set of discontinuity points of any totally monotone function on $I_a^b\subset\mathbb{R}^n$
lies on at most a countable set of hyperplanes of dimension $n-1$ parallel to the
coordinate axes. Given $i\in\{1,\dots,n\}$, denote by $Z_i$ the union of the set
of all rational points of the interval $[a_i,b_i]$, the two-point set $\{a_i,b_i\}$
and the set of those points $z_i\in[a_i,b_i]$, for which the hyperplane
  \begin{equation} \label{e:hyp}
H_i(z_i)=[a_1,b_1]\times\cdots\times[a_{i-1},b_{i-1}]\times\{z_i\}\times[a_{i+1},b_{i+1}]
\times\cdots\times[a_n,b_n]
  \end{equation}
contains points of discontinuity of $\nu$. Clearly, the sets $Z_i\subset[a_i,b_i]$
are countable and dense in $[a_i,b_i]$, and so, we may assume that
$Z_i=\{z_i(k)\}_{k=1}^\infty$.

\smallbreak
2. In order to apply the induction hypothesis, we need an estimate on the
$(n-1)$-dimensional total variation of any function $f=f_j$ from the sequence
$\{f_j\}$ `over the hyperplane' \eq{e:hyp} in the sense to be made precise
below. This is done as follows.

Let us fix $i\in\{1,\dots,n\}$ and set $1^i=(1,\dots,1,0,1,\dots,1)$, where
$0$ is the $i$-th coordinate of $1^i$ and the other coordinates of $1^i$ are
equal to~$1$. Note that $|1^i|=n-1$. Given $z_i\in Z_i$, we put
  \begin{equation} \label{e:bara}
\ov a\equiv\ov a(z_i)=(a_1,\dots,a_{i-1},z_i,a_{i+1},\dots,a_n).
  \end{equation}
The map $f_{1^i}^{\ov a}:I_a^b\lfloor1^i\to M$ with the base at $\ov a$,
truncated by $1^i$, is defined on the $(n-1)$-dimensional rectangle
$I_a^b\lfloor1^i\subset\mathbb{R}^{n-1}$ and given by: if $x\in I_a^b$,
then $x\lfloor1^i\in I_a^b\lfloor1^i$ and
  \begin{equation} \label{e:fova}
f_{1^i}^{\ov a}(x\lfloor1^i)=f(\ov a+1^i(x-\ov a))=f(x_1,\dots,x_{i-1},z_i,x_{i+1},\dots,x_n).
  \end{equation}
The $(n-1)$-dimensional total variation of $f_{1^i}^{\ov a}$ on $I_a^b\lfloor1^i$
is equal to
  \begin{equation} \label{e:tvn-1}
\mbox{TV}_{\!n-1}(f_{1^i}^{\ov a},I_a^b\lfloor1^i)=\sum_{0\ne\alpha\le1}
V_{|\alpha|}\bigl((f_{1^i}^{\ov a})_\alpha^{a\lfloor1^i},(I_a^b\lfloor1^i)\lfloor\alpha\bigr),
  \end{equation}
where the summation is taken over $(n-1)$-dimensional multiindices
$\alpha=(\alpha_1,\dots,\alpha_{n-1})$ s.t.\ $0_{n-1}\ne\alpha\le1_{n-1}$,
i.e., $\alpha\in\mathcal{A}(n-1)$ (this is the only instance and exception
when the summation is over $(n-1)$-dimensional multiindices).
Given $\alpha\in\mathcal{A}(n\!-\!1)$, we set
$\ov\alpha=(\alpha_1,\dots,\alpha_{i-1},0,\alpha_i,\dots,\alpha_{n-1})$, where
$0$ occupies the $i$-th place, and note that $\alpha=\ov\alpha\,\lfloor1^i$.
We have
  $$(f_{1^i}^{\ov a})_\alpha^{a\lfloor1^i}=f_{\ov\alpha}^{\ov a}\,\,\quad\mbox{on}\,\,\quad
    (I_a^b\lfloor1^i)\lfloor\alpha=I_a^b\lfloor\ov\alpha=I_{\ov a}^b\lfloor\ov\alpha.$$
In fact, given $x\in I_a^b$, we find $x\lfloor\ov\alpha=(x\lfloor1^i)\lfloor\alpha$ and
  \begin{eqnarray}
(f_{1^i}^{\ov a})_\alpha^{a\lfloor1^i}(x\lfloor\ov\alpha)&=&
  (f_{1^i}^{\ov a})_{\ov\alpha\lfloor1^i}^{a\lfloor1^i}\,((x\lfloor1^i)\lfloor\alpha)\nonumber\\[2pt]
&=&f_{1^i}^{\ov a}\bigl((a\lfloor1^i)+(\ov\alpha\,\lfloor1^i)[(x\lfloor1^i)-(a\lfloor1^i)]\bigr)
  \nonumber\\[2pt]
&=&f_{1^i}^{\ov a}\bigl([a+\ov\alpha\,(x-a)]\lfloor1^i\bigr)\nonumber\\[2pt]
&=&f\bigl(\ov a+1^i[a+\ov\alpha\,(x-a)-\ov a]\bigr).\label{e:xova}
  \end{eqnarray}
Since $\ov a+1^i(a-\ov a)=\ov a$ and $1^i\,\ov\alpha=\ov\alpha$, we get
  \begin{eqnarray*}
\ov a+1^i[a+\ov\alpha\,(x-a)-\ov a]&=&\ov a+1^i(a-\ov a)+1^i\,\ov\alpha\,(x-a)\\[2pt]
&=&\ov a+\ov\alpha\,(x-a)=\ov a+\ov\alpha\,(x-\ov a),
  \end{eqnarray*}
and so, the value \eq{e:xova} is equal to
  $$f(\ov a+\ov\alpha\,(x-\ov a))=f_{\ov\alpha}^{\ov a}(x\lfloor\ov\alpha).$$
It follows that the $|\alpha|$-th variation at the right-hand side of \eq{e:tvn-1}
is equal to
  $$V_{|\alpha|}\bigl((f_{1^i}^{\ov a})_\alpha^{a\lfloor1^i},(I_a^b\lfloor1^i)\lfloor\alpha\bigr)
    =V_{|\ov\alpha|}(f_{\ov\alpha}^{\ov a},I_{\ov a}^b\lfloor\ov\alpha).$$
Noting that the set $\mathcal{A}(n-1)$ is bijective to the set of those
$\ov\alpha\in\mathcal{A}(n)$, for which $0\ne\ov\alpha\le1^i$, and applying
Theorem~\ref{t:C} with $x=\ov a$, $y=b$ and $\gamma=1^i$, we get:
  \begin{eqnarray}
\mbox{TV}_{\!n-1}(f_{1^i}^{\ov a},I_a^b\lfloor1^i)&=&\sum_{0\ne\ov\alpha\le1^i}
  V_{|\ov\alpha|}(f_{\ov\alpha}^{\ov a},I_{\ov a}^b\lfloor\ov\alpha)
  =\mbox{TV}(f,I_{\ov a}^{\ov a+1^i(b-\ov a)})\nonumber\\[2pt]
&\le&\mbox{TV}(f,I_a^{\ov a+1^i(b-\ov a)})-\mbox{TV}(f,I_a^{\ov a})
  \le\mbox{TV}(f,I_a^b).\label{e:est}
  \end{eqnarray}
Thus, given $j\in\mathbb{N}$ and $i\in\{1,\dots,n\}$, setting back $f=f_j$, by virtue
of \eq{e:bara}, \eq{e:est} and \eq{e:ubv}, we find, for all $z_i\in Z_i$ and
$\ov a=\ov a(z_i)$:
  \begin{equation} \label{e:fjtv}
\mbox{TV}_{\!n-1}\bigl((f_j)_{1^i}^{\ov a(z_i)},I_a^b\lfloor1^i\bigr)\le C<\infty.
  \end{equation}

\smallbreak
3. Now, we make use of the diagonal processes. For $i=1$ and $z_1=z_i(1)=z_1(1)\in Z_1$
the sequence $\{(f_j)_{1^i}^{\ov a(z_i(1))}\}_{j=1}^\infty=\{(f_j)_{1^1}^{\ov a(z_1(1))}\}_{j=1}^\infty$
satisfies the uniform estimate \eq{e:fjtv} on the rectangle $I_a^b\lfloor1^1$ of dimension
$n-1$ and, since each map from this sequence is of the form \eq{e:fova} with
$z_i=z_1=z_1(1)$, then it follows from the assumptions of Theorem~\ref{t:first} that
the sequence under consideration is pointwise precompact on $I_a^b\lfloor1^1$.
By the induction hypothesis, the sequence $\{f_j\}$ contains a subsequence,
denoted by $\{f_j^1\}$, s.t.\ $(f_j^1)_{1^1}^{\ov a(z_1(1))}$ converges pointwise
on $I_a^b\lfloor1^1$ to a map from $I_a^b\lfloor1^1$ into $M$ of $(n-1)$-dimensional
finite total variation on $I_a^b\lfloor1^1$. Since, by \eq{e:fova},
  $$(f_j^1)_{1^1}^{\ov a(z_1(1))}(x_2,\dots,x_n)=(f_j^1)_{1^1}^{\ov a(z_1(1))}(x\lfloor1^1)
    =f_j^1(z_1(1),x_2,\dots,x_n)$$
with $x=(x_1,\dots,x_n)\in I_a^b$ and $x_i\in[a_i,b_i]$ for $i\in\{2,\dots,n\}$, then
the pointwise convergence above means, actually, that the sequence $\{f_j^1\}$ converges
pointwise on the hyperplane $H_1(z_1(1))=\{z_1(1)\}\times[a_2,b_2]\times\cdots\times[a_n,b_n]$.

Inductively, if $k\ge2$ and a subsequence $\{f_j^{k-1}\}_{j=1}^\infty$ of $\{f_j\}$,
which is pointwise convergent on $\bigcup_{l=1}^{k-1}H_1(z_1(l))$, is already chosen,
then the sequence $\{(f_j^{k-1})_{1^1}^{\ov a(z_1(k))}\}_{j=1}^\infty$ satisfies the
uniform estimate \eq{e:fjtv} on the rectangle $I_a^b\lfloor1^1$,
where $f_j$ is replaced by $f_j^{k-1}$ and $\ov a(z_i)$---by $\ov a(z_1(k))$.
Moreover, since, as above, the sequence is pointwise precompact on $I_a^b\lfloor1^1$,
then, by the induction hypothesis, there exists a subsequence $\{f_j^k\}_{j=1}^\infty$
of $\{f_j^{k-1}\}_{j=1}^\infty$ s.t.\ $(f_j^k)_{1^1}^{\ov a(z_1(k))}$ converges pointwise
on $I_a^b\lfloor1^1$ as $j\to\infty$ to a map from $I_a^b\lfloor1^1$ into $M$ of
$(n-1)$-dimensional finite total variation on $I_a^b\lfloor1^1$. Again, as above,
this pointwise convergence means that the sequence $\{f_j^k\}_{j=1}^\infty$ converges
pointwise on the hyperplane $H_1(z_1(k))$ and, as a consequence, on the set
$\bigcup_{l=1}^kH_1(z_1(l))$ as well. We infer that the diagonal sequence $\{f_j^j\}_{j=1}^\infty$,
which is a subsequence of the original sequence $\{f_j\}$, converges pointwise on
the set $H_1(Z_1)=\bigcup_{z_1\in Z_1}H_1(z_1)=\bigcup_{l=1}^\infty H_1(z_1(l))$;
in fact, given $(z_1,x_2,\dots,x_n)\in H_1(Z_1)$, we have $z_1=z_1(k)\in Z_1$ for
some $k\in\mathbb{N}$ and $(x_2,\dots,x_n)\in I_a^b\lfloor1^1$, and so, noting that
$\{f_j^j\}_{j=k}^\infty$ is a subsequence of $\{f_j^k\}_{j=1}^\infty$, we find that
  $$f_j^j(z_1,x_2,\dots,x_n)=(f_j^j)_{1^1}^{\ov a(z_1(k))}(x_2,\dots,x_n)$$
converges in $M$ as $j\to\infty$.

Let us denote the diagonal sequence $\{f_j^j\}_{j=1}^\infty$ extracted in the last
paragraph again by $\{f_j\}$. Then we let $i=2$, $z_2=z_i(1)=z_2(1)\in Z_2$ and,
beginning with the sequence
$\{(f_j)_{1^i}^{\ov a(z_i(1))}\}_{j=1}^\infty=\{(f_j)_{1^2}^{\ov a(z_2(1))}\}_{j=1}^\infty$,
apply the above arguments of this step. Doing this, we will end up with a diagonal
sequence, a subsequence of the original sequence $\{f_j\}$, again denoted by $\{f_j\}$,
which converges pointwise on $H_1(Z_1)\cup H_2(Z_2)$. Now suppose that for some
$i\in\{2,\dots,n-1\}$ we have already extracted a (diagonal) subsequence of $\{f_j\}$,
again denoted by $\{f_j\}$, which converges pointwise on the set
$H_1(Z_1)\cup\dots\cup H_{i-1}(Z_{i-1})$. Then we let $z_i=z_i(1)\in Z_i$ and apply
the above arguments of this step to the sequence $\{(f_j)_{1^i}^{\ov a(z_i(1))}\}_{j=1}^\infty$:
a subsequence of the original sequence $\{f_j\}$ converges pointwise on the set
$H_1(Z_1)\cup\dots\cup H_i(Z_i)$. In this way after finitely many steps we obtain
a subsequence of the original sequence $\{f_j\}$, again denoted by $\{f_j\}$,
which converges pointwise on the set $H(Z)=\bigcup_{i=1}^nH_i(Z_i)$.

\smallbreak
4. Finally, let us show that the sequence $\{f_j\}$ from the end of Step~3 converges
at each point $y\in I_a^b\setminus H(Z)$. Note that $y$ is a point of continuity of
the function $\nu$ from \eq{e:nujnu} s.t.\ its coordinates $a_i<y_i<b_i$ are
irrational for all $i\in\{1,\dots,n\}$. Due to the density of $H(Z)$ in $I_a^b$,
the continuity of $\nu$ at $y$ and properties of totally monotone functions,
given $\varepsilon>0$, there exists $x=x(\varepsilon)\in H(Z)$ with $x<y$ s.t.\
$0\le\nu(y)-\nu(x)\le\varepsilon$. By virtue of \eq{e:nujnu}, choose a number
$j_0(\varepsilon)\in\mathbb{N}$ s.t.\ $|\nu_j(y)-\nu(y)|\le\varepsilon$ and
$|\nu(x)-\nu_j(x)|\le\varepsilon$ for all $j\ge j_0(\varepsilon)$.
By Theorems \ref{t:B} and \ref{t:C} with $\gamma=1$, for all $j\ge j_0(\varepsilon)$
we have:
  \begin{eqnarray*}
d(f_j(x),f_j(y))&\le&\mbox{TV}(f_j,I_x^y)\le\mbox{TV}(f_j,I_a^y)-\mbox{TV}(f_j,I_a^x)
  =\nu_j(y)-\nu_j(x)\\[2pt]
&\le&|\nu_j(y)-\nu(y)|+(\nu(y)-\nu(x))+|\nu(x)-\nu_j(x)|\le3\varepsilon.
  \end{eqnarray*}
Since $x\in H(Z)$ and, as it was shown in Step~3, the sequence $\{f_j(x)\}_{j=1}^\infty$
is convergent in $M$, it is Cauchy, and so, there exists a number $j_1(\varepsilon)\in\mathbb{N}$
s.t.\ $d(f_j(x),f_{j'}(x))\le\varepsilon$ for all $j\ge j_1(\varepsilon)$ and
$j'\ge j_1(\varepsilon)$. It follows that if
$J(\varepsilon)=\max\{j_0(\varepsilon),j_1(\varepsilon)\}$, $j\ge J(\varepsilon)$ and
$j'\ge J(\varepsilon)$, then we have:
  \begin{eqnarray*}
d(f_j(y),f_{j'}(y))&\le&d(f_j(y),f_j(x))+d(f_j(x),f_{j'}(x))+d(f_{j'}(x),f_{j'}(y))\\[2pt]
&\le&3\varepsilon+\varepsilon+3\varepsilon=7\varepsilon.
  \end{eqnarray*}
Thus, the sequence $\{f_j(y)\}_{j=1}^\infty$ is Cauchy in the metric space $M$, and so,
since it is also precompact by the assumption, it is convergent in~$M$.

It follows from here and the end of Step~3 that the sequence $\{f_j(y)\}_{j=1}^\infty$
converges in $M$ at each point $y\in(I_a^b\setminus H(Z))\cup H(Z)=I_a^b$, i.e., the
sequence $\{f_j\}$, which is a subsequence of the original sequence $\{f_j\}$,
converges pointwise on $I_a^b$. Let us denote the pointwise limit of $\{f_j\}$ by
$f:I_a^b\to M$. Then, by virtue of Theorem~\ref{t:E} and assumption \eq{e:ubv}, we find
  $$\mbox{TV}(f,I_a^b)\le\liminf_{j\to\infty}\mbox{TV}(f_j,I_a^b)\le C,$$
and so, $f\in\mbox{BV}(I_a^b;M)$.

This completes the proof of Theorem~\ref{t:first}.
\qed\end{pfirst}

\begin{rem} \label{r:cb}
In Theorem~\ref{t:first} the precompactness of the sets $\{f_j(x)\}_{j=1}^\infty$
at all points $x\in I_a^b$ cannot be replaced by the closedness and boundedness
even at a single point of $I_a^b$. The corresponding examples for maps of one
variable are constructed in \cite[Section~3]{jmaa00}, \cite[Section~5]{msb98}
and \cite[Section~1]{Sovae} and can be easily adapted for maps of several variables.
\end{rem}

\begin{pweak}
The proof is adapted for the situation under consideration from the proof of
Theorem~7 from \cite{jmaa05}.

1. In this step we show that, given $j\in\mathbb{N}$ and $u^*\in M^*$, we have:
  \begin{equation} \label{e:W.3}
\mbox{TV}(\langle f_j(\cdot),u^*\rangle,I_a^b)\le\mbox{TV}(f_j,I_a^b)\|u^*\|\le C\|u^*\|,
  \end{equation}
where the function $\langle f_j(\cdot),u^*\rangle:I_a^b\to\mathbb{K}$ is given by
$\langle f_j(\cdot),u^*\rangle(x)=\langle f_j(x),u^*\rangle$, $x\in I_a^b$, and $C$
is the constant from \eq{e:ubv}.

In fact, given $0\ne\alpha\le1$ and $x,y\in I_a^b$ with $x<y$, by virtue of \eq{e:2.4*}
where $d(u,v)$ is replaced by the absolute value $|u-v|$ in $\mathbb{K}$ and later on---%
by the norm in $M$, we get:
  \begin{eqnarray*}
\mbox{md}_{|\alpha|}(\langle f_j(\cdot),u^*\rangle_\alpha^a,I_x^y\lfloor\alpha)
  &=&\biggl|\sum_{0\le\theta\le\alpha}(-1)^{|\theta|}
  \langle f_j(a\!+\!\alpha\,(x\!-\!a)\!+\!\theta(y\!-\!x)),u^*\rangle\biggr|\\
&=&\biggl|\Bigl\langle\sum_{0\le\theta\le\alpha}(-1)^{|\theta|}
  f_j(a\!+\!\alpha\,(x\!-\!a)\!+\!\theta(y\!-\!x)),u^*\Bigr\rangle\biggr|\\
&\le&\biggl\|\sum_{0\le\theta\le\alpha}(-1)^{|\theta|}
  f_j(a\!+\!\alpha\,(x\!-\!a)\!+\!\theta(y\!-\!x))\biggr\|\cdot\|u^*\|\\[3pt]
&=&\mbox{md}_{|\alpha|}((f_j)_\alpha^a,I_x^y\lfloor\alpha)\|u^*\|.
  \end{eqnarray*}
It follows that if $\mathcal{P}=\{x[\sigma]\}_{\sigma=0}^\kappa$ is a net partition
of $I_a^b$, then $\mathcal{P}\lfloor\alpha=\{x[\sigma]\lfloor\alpha\}_{\sigma\lfloor\alpha=0}%
^{\kappa\lfloor\alpha}$ is a net partition of $I_a^b\lfloor\alpha$, and so, setting
$x=x[\sigma-1]$ and $y=x[\sigma]$ in the calculations above, we find
  \begin{eqnarray*}
\sum_{1\lfloor\alpha\le\sigma\lfloor\alpha\le\kappa\lfloor\alpha}
  \!\!\!\!\!\!\!\!\mbox{md}_{|\alpha|}\bigl(\langle f_j(\cdot),u^*\rangle_\alpha^a,
  I_{x[\sigma-1]}^{x[\sigma]}\lfloor\alpha\bigr)
&\le&\sum_{1\lfloor\alpha\le\sigma\lfloor\alpha\le\kappa\lfloor\alpha}
  \!\!\!\!\!\!\!\!\mbox{md}_{|\alpha|}\bigl((f_j)_\alpha^a,
  I_{x[\sigma-1]}^{x[\sigma]}\lfloor\alpha\bigr)\|u^*\|\\
&\le&V_{|\alpha|}((f_j)_\alpha^a,I_a^b\lfloor\alpha)\|u^*\|,
  \end{eqnarray*}
the summation over $\sigma\lfloor\alpha$ being taken only over those coordinates $\sigma_i$
in the range $1\le\sigma_i\le\kappa_i$ with $i\in\{1,\dots,n\}$, for which $\alpha_i=1$.
Since $\mathcal{P}$ is an arbitrary partition of $I_a^b$, we get:
  $$V_{|\alpha|}\bigl(\langle f_j(\cdot),u^*\rangle_\alpha^a,I_a^b\lfloor\alpha\bigr)\le
    V_{|\alpha|}\bigl((f_j)_\alpha^a,I_a^b\lfloor\alpha\bigr)\|u^*\|,$$
and so, inequality \eq{e:W.3} follows from the definition of the total variation.

Moreover, by virtue of \eq{e:pou}, we have:
  \begin{equation} \label{e:ufjx}
|\langle f_j(x),u^*\rangle|\le\|f_j(x)\|\cdot\|u^*\|\le c(x)\|u^*\|,\qquad
x\in I_a^b,\quad u^*\in M^*,
  \end{equation}
and so, the sequence $\{\langle f_j(\cdot),u^*\rangle\}_{j=1}^\infty$ of functions from
$I_a^b$ into (metric semigroup) $\mathbb{K}$ is pointwise bounded on $I_a^b$ and, hence,
pointwise precompact for each $u^*\in M^*$.

Taking this and \eq{e:W.3} into account and applying Theorem~\ref{t:first} to the
sequence $\{\langle f_j(\cdot),u^*\rangle\}_{j=1}^\infty$ for any given $u^*\in M^*$,
we extract a subsequence of $\{f_j\}$, denoted by $\{f_{j,u^*}\}$ (which depends on $u^*$
in general), and find a function $g_{u^*}\in\mbox{BV}(I_a^b;\mathbb{K})$ satisfying
$\mbox{TV}(g_{u^*},I_a^b)\le C\|u^*\|$ s.t.\ $\langle f_{j,u^*}(x),u^*\rangle\to g_{u^*}(x)$
in $\mathbb{K}$ as $j\to\infty$ for all $x\in I_a^b$.

\smallbreak
2. Making use of the diagonal process and the separability of $M^*$, let us get
rid of the dependence of $\{f_{j,u^*}\}$ on $u^*\in M^*$. Let $\{u_k^*\}_{k=1}^\infty$
be a countable dense subset of $M^*$. By Step~1, for $u^*=u^*_1$ we get a subsequence
$\{f_j^{(1)}\}=\{f_{j,u_1^*}\}_{j=1}^\infty$ of the original sequence $\{f_j\}$ and
a function $g_{u_1^*}\in\mbox{BV}(I_a^b;\mathbb{K})$ satisfying
$\mbox{TV}(g_{u_1^*},I_a^b)\le C\|u_1^*\|$ s.t.\
$\langle f_j^{(1)}(x),u_1^*\rangle\to g_{u_1^*}(x)$ in $\mathbb{K}$ for all $x\in I_a^b$.
Inductively, if $k\ge2$ and a subsequence $\{f_j^{(k-1)}\}_{j=1}^\infty$ of $\{f_j\}$
is already chosen, then by virtue of \eq{e:W.3} and \eq{e:ufjx}, we have:
  $$\mbox{TV}(\langle f_j^{(k-1)}(\cdot),u_k^*\rangle,I_a^b)\le C\|u_k^*\|$$
and $|\langle f_j^{(k-1)}(x),u_k^*\rangle|\le c(x)\|u_k^*\|$, $x\in I_a^b$, for all
$j\in\mathbb{N}$. By Theorem~\ref{t:first}, applied to the sequence
$\{\langle f_j^{(k-1)}(\cdot),u_k^*\rangle\}_{j=1}^\infty$, there exist a subsequence
$\{f_j^{(k)}\}_{j=1}^\infty$ of $\{f_j^{(k-1)}\}_{j=1}^\infty$ and a function
$g_{u_k^*}\in\mbox{BV}(I_a^b;\mathbb{K})$ satisfying
$\mbox{TV}(g_{u_k^*},I_a^b)\le C\|u_k^*\|$ s.t.\
$\langle f_j^{(k)}(x),u_k^*\rangle\to g_{u_k^*}(x)$ in $\mathbb{K}$ as $j\to\infty$
for all $x\in I_a^b$. Then the diagonal sequence $\{f_j^{(j)}\}_{j=1}^\infty$,
again denoted by $\{f_j\}$, is a subsequence of the original sequence and satisfies
the condition:
  \begin{equation} \label{e:W5}
\langle f_j(x),u_k^*\rangle\to g_{u_k^*}(x)\quad\mbox{as $j\to\infty$ for all $x\in I_a^b$
and $k\in\mathbb{N}$.}
  \end{equation}

\smallbreak
3. Now, given $u^*\in M^*$ and $x\in I_a^b$, let us show that the sequence
$\{\langle f_j(x),u^*\rangle\}_{j=1}^\infty$ is Cauchy in $\mathbb{K}$. Taking
into account \eq{e:W5} we may assume that $u^*\ne u_k^*$ for all $k\in\mathbb{N}$.
Let $\varepsilon>0$ be arbitrary. By the density of $\{u_k^*\}_{k=1}^\infty$ in $M^*$,
there exists $k=k(\varepsilon)\in\mathbb{N}$ s.t.\
$\|u^*-u_k^*\|\le\varepsilon/(1+4c(x))$. By \eq{e:W5}, there exists
$j_0=j_0(\varepsilon)\in\mathbb{N}$ s.t.\
$|\langle f_j(x),u_k^*\rangle-\langle f_{j'}(x),u_k^*\rangle|\le\varepsilon/2$
for all $j\ge j_0$ and $j'\ge j_0$. It follows that for such $j$ and $j'$ we have:
  \begin{eqnarray*}
|\langle f_j(x),u^*\rangle-\langle f_{j'}(x),u^*\rangle|&\le&
  |\langle f_j(x)-f_{j'}(x),u^*-u_k^*\rangle|\\[2pt]
&&\quad\,+|\langle f_j(x),u_k^*\rangle-\langle f_{j'}(x),u_k^*\rangle|\\[2pt]
&\le&\|f_j(x)-f_{j'}(x)\|\cdot\|u^*-u_k^*\|+\displaystyle\frac\varepsilon2\\[2pt]
&\le&\displaystyle 2c(x)\frac{\varepsilon}{1+4c(x)}+\frac\varepsilon2\le\varepsilon.
  \end{eqnarray*}
Thus, $\{\langle f_j(x),u^*\rangle\}_{j=1}^\infty$ is Cauchy in $\mathbb{K}$ and, hence,
there exists an element of $\mathbb{K}$ denoted by $g_{u^*}(x)$ s.t.\
$\langle f_j(x),u^*\rangle\to g_{u^*}(x)$ in $\mathbb{K}$ as $j\to\infty$.
In other words, we have shown that for each $u^*\in M^*$ there exists a function
$g_{u^*}:I_a^b\to\mathbb{K}$ satisfying (cf.\ Theorem~\ref{t:E} and \eq{e:W.3})
  $$\mbox{TV}(g_{u^*},I_a^b)\le\liminf_{j\to\infty}\mbox{TV}(\langle f_j(\cdot),u^*\rangle,I_a^b)
    \le C\|u^*\|$$
(and so, $g_{u^*}\in\mbox{BV}(I_a^b;\mathbb{K})$) and
  \begin{equation} \label{e:cuec}
\lim_{j\to\infty}\langle f_j(x),u^*\rangle=g_{u^*}(x)\quad\mbox{in $\mathbb{K}$ for all
$x\in I_a^b$ and $u^*\in M^*$.}
  \end{equation}

\smallbreak
4. Let us prove \eq{e:wec}, i.e., $f_j(x)$ converges weakly in $M$ as $j\to\infty$
for all $x\in I_a^b$. By the reflexivity of $M$, we have
$f_j(x)\in M=M^{**}\equiv\mbox{L}(M^*;\mathbb{K})$ for all $j\in\mathbb{N}$.
Define the functional $G_x:M^*\to\mathbb{K}$ by $G_x(u^*)=g_{u^*}(x)$ for all
$u^*\in M^*$. By virtue of \eq{e:cuec}, we get
  $$\lim_{j\to\infty}\langle f_j(x),u^*\rangle=g_{u^*}(x)=G_x(u^*)\quad
    \mbox{for \,\,all}\quad u^*\in M^*,$$
i.e., the sequence $\{f_j(x)\}_{j=1}^\infty\subset\mbox{L}(M^*;\mathbb{K})$ converges
pointwise on $M^*$ to the operator $G_x:M^*\to\mathbb{K}$. By the Banach-Steinhaus
uniform boundedness principle, $G_x\in\mbox{L}(M^*;\mathbb{K})=M$ and
$\|G_x\|\le\liminf_{j\to\infty}\|f_j(x)\|$. Setting $f(x)=G_x$ for all $x\in I_a^b$,
we find that $f:I_a^b\to M$ and
  \begin{equation} \label{e:W-7}
\lim_{j\to\infty}\langle f_j(x),u^*\rangle=G_x(u^*)=\langle G_x,u^*\rangle=
\langle f(x),u^*\rangle\quad\mbox{in}\quad\mathbb{K}
  \end{equation}
for all $u^*\in M^*$ and $x\in I_a^b$, and so, $f_j(x)\wto f(x)$ in $M$ as $j\to\infty$
for all $x\in I_a^b$, which proves \eq{e:wec}.

\smallbreak
5. It remains to show that $f\in\mbox{BV}(I_a^b;M)$ and $\mbox{TV}(f,I_a^b)\le C$.
By \eq{e:W-7}, we have: if $x,y\in I_a^b$ with $x<y$ and $0\ne\alpha\le1$, then
  $$\sum_{0\le\theta\le\alpha}\!(-1)^{|\theta|}f_j\bigl(
    a\!+\!\alpha\,(x\!-\!a)\!+\!\theta(y\!-\!x)\bigr)\wto
    \sum_{0\le\theta\le\alpha}\!(-1)^{|\theta|}f\bigl(
    a\!+\!\alpha\,(x\!-\!a)\!+\!\theta(y\!-\!x)\bigr)$$
in $M$ as $j\to\infty$, and so, by virtue of \eq{e:2.4*} and the remarks preceding
Theorem~\ref{t:weak},
  \begin{eqnarray}
\mbox{md}_{|\alpha|}(f_\alpha^a,I_x^y\lfloor\alpha)&=&\biggl\|\sum_{0\le\theta\le\alpha}
  (-1)^{|\theta|}f\bigl(a\!+\!\alpha\,(x\!-\!a)\!+\!\theta(y\!-\!x)\bigr)\biggr\|\nonumber\\
&\le&\liminf_{j\to\infty}\,\biggl\|\sum_{0\le\theta\le\alpha}
  (-1)^{|\theta|}f_j\bigl(a\!+\!\alpha\,(x\!-\!a)\!+\!\theta(y\!-\!x)\bigr)\biggr\|\nonumber\\[3pt]
&=&\liminf_{j\to\infty}\,\mbox{md}_{|\alpha|}((f_j)_\alpha^a,I_x^y\lfloor\alpha).\label{e:Wei}
  \end{eqnarray}
Arguing as in Step~2 of the proof of Theorem~\ref{t:E}, making use of the inequality
\eq{e:Wei}, which coincides with \eq{e:almin} (see p.~\pageref{pg:li}), and taking into
account \eq{e:ubv}, we get:
  $$\mbox{TV}(f,I_a^b)\le\liminf_{j\to\infty}\,\mbox{TV}(f_j,I_a^b)\le C.$$
\par This completes the proof of Theorem~\ref{t:weak}.
\qed\end{pweak}

\begin{rem} \label{r:ca}
Note that instead of condition \eq{e:pou} in Theorem~\ref{t:weak} we may assume only that
the value $c(a)=\sup_{j\in\mathbb{N}}\|f_j(a)\|$ is finite. In fact, it follows from
Theorem~\ref{t:B} and condition \eq{e:ubv} that, given $x\in I_a^b$ and $j\in\mathbb{N}$,
  $$\|f_j(x)\|\le\|f_j(a)\|+\|f_j(x)-f_j(a)\|\le c(a)+\mbox{TV}(f_j,I_a^x)\le c(a)+C.$$
\end{rem}

\section{Proof of Theorem \ref{t:B}} \label{s:prB}

In order to prove Theorem~\ref{t:B}, we need three more Lemmas \ref{l:binco},
\ref{l:metsem} and \ref{l:funcah}. The following equality will be needed in
Lemma~\ref{l:binco} (cf.\ \cite[Part~I, equality (3.4)]{na05}): given two
multiindices $0\le\beta\le\gamma\le1$, we have:
  \begin{equation} \label{e:d.1}
|\{\alpha:\mbox{$\beta\le\alpha\le\gamma$ and $|\alpha|=i$}\}|=C_{|\gamma|-|\beta|}^{\,i-|\beta|}
\quad\,\,\mbox{for \,\,all}\quad\,\,|\beta|\le i\le|\gamma|,
  \end{equation}
where $|A|$ denotes the number of elements in the set $A$ and, given $0\le j\le m$,
$C_m^{\,j}=\binom{m}{j}=\frac{m!}{j!(m-j)!}$ is the usual binomial coefficient (with $0!=1$).
Also, recall (cf.\ Section~\ref{s:defs}) that a multiindex $\alpha$ is said to be even
(odd) if $\alpha\in\mathcal{E}(n)$ ($\alpha\in\mathcal{O}(n)$, respectively).

\begin{lem} \label{l:binco}
{\rm(a)} Given $m\in\mathbb{N}$ and integer $0\le k\le m-1$, we have\/{\rm:}
  $$\sum_{i\ge k/2}^{\le m/2}C_{m-k}^{\,2i-k}\,\,=\sum_{i\ge (k+1)/2}^{\,\le(m+1)/2}\!\!
    C_{m-k}^{\,2i-1-k}\,=2^{m-k-1},$$
where the summations are taken over {\sl integer} $i$ in the ranges specified.

{\rm(b)} Given two multiindices $0\le\beta\le\gamma\le1$ with $\beta\ne\gamma$, we have\/{\rm:}
  $$|\{\mbox{\rm even $\alpha$}:\beta\le\alpha\le\gamma\}|
    =|\{\mbox{\rm odd $\alpha$}:\beta\le\alpha\le\gamma\}|.$$
\end{lem}

\begin{pf}
(a) Since the case $m=1$ is clear, we suppose that $m\ge2$. By the binomial formula,
$0=(1-1)^{m-k}=\sum_{j=0}^{m-k}(-1)^jC_{m-k}^{\,j}$, which is equal to
  $$\sum_{j=0}^{(m-k-1)/2}\!\!C_{m-k}^{\,2j}-\sum_{j=1}^{(m-k+1)/2}\!\!C_{m-k}^{\,2j-1}\quad\,\,
    \mbox{if $m$ and $k$ are of different evenness,}$$
and
  $$\sum_{j=0}^{(m-k)/2}C_{m-k}^{\,2j}-\sum_{j=1}^{(m-k)/2}C_{m-k}^{\,2j-1}\quad\,\,
    \mbox{if $m$ and $k$ are of the same evenness.}$$
If $k$ is even, we change the summation index $j\mapsto i=(2j+k)/2$ in these sums
and get, respectively:
  $$0=\sum_{i=k/2}^{(m-1)/2}C_{m-k}^{\,2i-k}-\sum_{i=(k+2)/2}^{(m+1)/2}C_{m-k}^{\,2i-1-k}
    \qquad\mbox{if $m$ is odd,}$$
  $$0=\sum_{i=k/2}^{m/2}C_{m-k}^{\,2i-k}-\sum_{i=(k+2)/2}^{m/2}C_{m-k}^{\,2i-1-k}
    \qquad\mbox{if $m$ is even,}$$
while if $k$ is odd, we change the summation index $j\mapsto i=(2j+1+k)/2$ in the
first sum and $j\mapsto i=(2j-1+k)/2$ in the second sum and get, respectively:
  $$0=\sum_{i=(k+1)/2}^{m/2}C_{m-k}^{\,2i-1-k}-\sum_{i=(k+1)/2}^{m/2}C_{m-k}^{\,2i-k}
    \qquad\mbox{if $m$ is even,}$$
  $$0=\sum_{i=(k+1)/2}^{(m+1)/2}C_{m-k}^{\,2i-1-k}-\sum_{i=(k+1)/2}^{(m-1)/2}C_{m-k}^{\,2i-k}
    \qquad\mbox{if $m$ is odd.}$$
The second equality in (a) follows from the equality $(1+1)^{m-k}=\sum_{j=0}^{m-k}C_{m-k}^{\,j}$.

\smallbreak
{\it Remark}. The first equality in Lemma~\ref{l:binco}(a) can be written also as
  $$\sum_{i=[(k+2)/2]}^{[(m+1)/2]}C_{m-k}^{\,2i-1-k}=\sum_{i=[(k+1)/2]}^{[m/2]}C_{m-k}^{\,2i-k},
    \quad\,\,m\in\mathbb{N},\quad 0\le k\le m-1,$$
where $[r]$ designates the integer part of $r\in\mathbb{R}$. However, we prefer the
equality in Lemma~\ref{l:binco}(a) since it is more simple and suggestive.

\smallbreak
(b) By virtue of \eq{e:d.1}, the left-hand side of the equality is equal to
  $$\bigl|\{\alpha:\mbox{$\beta\le\alpha\le\gamma$ and $|\alpha|=2i$ for all $i$ s.t.\
    $|\beta|\le2i\le|\gamma|$}\}\bigr|
   =\displaystyle\sum_{i\ge|\beta|/2}^{\le|\gamma|/2}C_{|\gamma|-|\beta|}^{\,2i-|\beta|},$$
and the right-hand side of the equality is equal to
  \begin{eqnarray*}
&\bigl|\{\alpha:\mbox{$\beta\le\alpha\le\gamma$ and $|\alpha|=2i-1$ for all $i$ s.t.\
  $|\beta|\le2i-1\le|\gamma|$}\}\bigr|&\\
&=\displaystyle\sum_{i\ge(|\beta|+1)/2}^{\le(|\gamma|+1)/2}C_{|\gamma|-|\beta|}^{\,2i-1-|\beta|}.&
  \end{eqnarray*}
It remains to put $m=|\gamma|$ and $k=|\beta|$, note that $k<m$ and apply the equality
from the previous assertion (a).
\qed\end{pf}

If $(M,d,+)$ is a metric semigroup, then, by virtue of the triangle inequality for $d$
and the translation invariance of metric $d$ on $M$, we have, for all $u,v,u',v'\in M$:
  \begin{eqnarray}
d(u,v)&\le&d(u',v')+d(u+u',v+v'),\nonumber\\[2pt]
d(u+u',v+v')&\le&d(u,v)+d(u',v').\label{e:ms.2}
  \end{eqnarray}
Inequality \eq{e:ms.2} yields that the addition operation $(u,v)\mapsto u+v$ is
a continuous map from $M\times M$ into $M$. More generally, if $u_j\to u$,
$v_j\to v$, $u_j'\to u'$ and $v_j'\to v'$ as $j\to\infty$ (convergence of sequences in $M$),
then $\lim_{j\to\infty}d(u_j+v_j,u_j'+v_j')=d(u+v,u'+v')$.

\begin{lem} \label{l:metsem}
If $m\in\mathbb{N}$, $u,v\in M$, $\{u_j\}_{j=1}^m$, $\{v_j\}_{j=1}^m\subset M$ and
  \begin{equation} \label{e:4s}
\sum_{i=1}^{\le m/2}u_{2i}\,+\,u\,+\sum_{i=1}^{\le(m+1)/2}\!\!v_{2i-1}=
\sum_{i=1}^{\le m/2}v_{2i}\,+\,v\,+\sum_{i=1}^{\le(m+1)/2}\!\!u_{2i-1},
  \end{equation}
then
  \begin{equation} \label{e:4d}
d(u,v)\le\sum_{j=1}^md(u_j,v_j).
  \end{equation}
\end{lem}

\begin{pf}
Observe that if $u+\ell_1+\cdots+\ell_k=v+r_1+\cdots+r_k$ for some $k\in\mathbb{N}$ and
$\{\ell_i,r_i\}_{i=1}^k\subset M$, then $d(u,v)\le\sum_{i=1}^kd(\ell_i,r_i)$. In fact,
by the translation invariance of $d$ and inequality \eq{e:ms.2}, we have:
  \begin{eqnarray}
d(u,v)&=&d\biggl(u+\sum_{i=1}^k\ell_i,v+\sum_{i=1}^k\ell_i\biggr)=
  d\biggl(v+\sum_{i=1}^kr_i,v+\sum_{i=1}^k\ell_i\biggr)\nonumber\\
&=&d\biggl(\sum_{i=1}^kr_i,\sum_{i=1}^k\ell_i\biggr)\le\sum_{i=1}^kd(r_i,\ell_i).\nonumber
  \end{eqnarray}
Applying this observation and equality \eq{e:4s}, we get:
  $$d(u,v)\le\sum_{i=1}^{\le m/2}d(u_{2i},v_{2i})+\sum_{i=1}^{\le(m+1)/2}d(v_{2i-1},u_{2i-1})
    =\sum_{j=1}^md(u_j,v_j).\qquad\square$$
\end{pf}

\begin{rem} \label{r:u+v}
In particular, {\rm(}in\/{\rm)}equalities\/ \eq{e:4s} and\/ \eq{e:4d} hold for {\sl odd\/} $m$ if
  \begin{equation} \label{e:odm}
u+\sum_{i=1}^{(m-1)/2}u_{2i}=\sum_{i=1}^{(m+1)/2}u_{2i-1}\quad\mbox{and}\quad
v+\sum_{i=1}^{(m-1)/2}v_{2i}=\sum_{i=1}^{(m+1)/2}v_{2i-1},
  \end{equation}
and for {\sl even\/} $m$ if either
  \begin{equation} \label{e:evm1}
u+\sum_{i=1}^{m/2}u_{2i}=v+\sum_{i=1}^{m/2}u_{2i-1}\quad\mbox{and}\quad
\sum_{i=1}^{m/2}v_{2i}=\sum_{i=1}^{m/2}v_{2i-1},
  \end{equation}
or
  \begin{equation} \label{e:evm2}
\sum_{i=1}^{m/2}u_{2i}=v+\sum_{i=1}^{m/2}u_{2i-1}\quad\mbox{and}\quad
\sum_{i=1}^{m/2}v_{2i}=u+\sum_{i=1}^{m/2}v_{2i-1}.
  \end{equation}
\end{rem}

In the next lemma we set $\mathcal{A}_0\equiv\mathcal{A}_0(n)=\{\theta\in\mathbb{N}_0^n:\theta\le1\}$.
Also, we stick to the following conventions: `$u\doteq0$' will mean that $u$ is omitted in the
formula under consideration (especially in a metric semigroup with no zero), and a sum over
the empty set is also omitted in any context (i.e., $\sum_{\varnothing}\doteq0$).

\begin{lem} \label{l:funcah}
Given a map $h:\mathcal{A}_0\to M$ and a multiindex $\gamma\in\mathcal{A}_0$, we have\/{\rm:}
  \begin{equation} \label{e:d.7.1}
\sum_{\even\alpha\le\gamma}\,\sum_{\even\theta\le\alpha}h(\theta)
=c_\gamma+\sum_{\odd\alpha\le\gamma}\,\sum_{\even\theta\le\alpha}h(\theta),
  \end{equation}
where $c_\gamma\doteq0$ if $\gamma$ is odd, and $c_\gamma=h(\gamma)$ if $\gamma$ is even, and
  \begin{equation} \label{e:d.7.2}
\sum_{\odd\alpha\le\gamma}\,\sum_{\odd\theta\le\alpha}h(\theta)
=d_\gamma+\sum_{\even\alpha\le\gamma}\,\sum_{\odd\theta\le\alpha}h(\theta),
  \end{equation}
where $d_\gamma=h(\gamma)$ if $\gamma$ is odd, and $d_\gamma\doteq0$ if $\gamma$ is even.
\end{lem}

\begin{pf}
0. Denote by $\mathcal{L}$ (by $\mathcal{R}$) the set of all `admissible' $\theta$'s at the
left (right) hand side of the equality under consideration and, given $\theta\in\mathcal{L}$
(and $\theta\in\mathcal{R}$), by $L(\theta)$ (and by $R(\theta)$)---the multiplicity of
the term $h(\theta)$ at the left (and right) hand sum(s). Then the equality can be
rewritten as
  \begin{equation} \label{e:LRh}
\sum_{\theta\in\mathcal{L}}L(\theta)h(\theta)=\sum_{\theta\in\mathcal{R}}R(\theta)h(\theta),
  \end{equation}
where $L(\theta)h(\theta)$ denotes the sum of terms of the form $h(\theta)$ taken $L(\theta)$
times (and likewise for $R(\theta)h(\theta)$). In what follows in order to prove \eq{e:LRh},
we show that $\mathcal{L}=\mathcal{R}$ and $L(\theta)=R(\theta)$ for all
$\theta\in\mathcal{L}=\mathcal{R}$.

We divide the proof into four steps for clarity.

In the first two steps we let $\gamma$ be odd (i.e., $0\le\gamma\le1$ and $|\gamma|$ is odd).

1. Let us establish \eq{e:d.7.1}. We have $\LC=\{\mbox{even $\theta$}:\exists\,
\mbox{even $\alpha\le\gamma$ s.t.\ $\theta\le\alpha$}\}$, i.e.,
$\LC=\{\mbox{even }\theta:\theta\le\gamma\}$, and
$\RC=\{\mbox{even }\theta:\exists\,\mbox{odd }\alpha\le\gamma\,\mbox{s.t.\,}\theta\le\alpha\}$.
The sets $\LC$ and $\RC$ are nonempty ($0\in\LC$ and $0\in\RC$) and $\LC=\RC$. In fact,
the inclusion $\LC\supset\RC$ is clear, and so, we let $\theta\in\LC$. Since $\theta$
is even, $\gamma$ is odd and $\theta\le\gamma$, there exists $i\in\{1,\dots,n\}$ s.t.\
$\theta_i=0$ and $\gamma_i=1$. We set $\alpha=(\theta_1,\dots,\theta_{i-1},1,\theta_{i+1},\dots,\theta_n)$.
It follows that $\alpha\le\gamma$, $|\alpha|=|\theta|+1$ is odd and $\theta\le\alpha$,
and so, $\theta\in\RC$.

Given $\theta\in\LC=\RC$, we find $\theta\ne\gamma$,
$L(\theta)=|\{\mbox{even }\alpha:\theta\le\alpha\le\gamma\}|$ and
$R(\theta)=|\{\mbox{odd }\alpha:\theta\le\alpha\le\gamma\}|$. By Lemma~\ref{l:binco}(b),
$L(\theta)=R(\theta)$, and so, \eq{e:LRh} holds implying \eq{e:d.7.1} with $c_\gamma\doteq0$.

\smallbreak
2. Let us prove \eq{e:d.7.2}. If $|\gamma|=1$, then the equality is immediate: the left-hand
side is equal to $h(\gamma)=d_\gamma$, while the double sum at the right is omitted
(in fact, even $\alpha\le\gamma$ implies $\alpha=0$, and so, no odd $\theta$ s.t.\
$\theta\le0$ exists). Now, if $|\gamma|>1$, then $\LC=\{\mbox{odd }\theta:\theta\le\gamma\}$
and $\RC=\{\mbox{odd }\theta:\exists\,\mbox{even $\alpha\le\gamma$}\,\mbox{s.t.\,$\theta\le\alpha$}\}%
\cup\{\gamma\}$ (disjoint union), and $\LC=\RC$. Let $\theta\in\LC=\RC$. If $\theta\ne\gamma$, then
$L(\theta)=|\{\mbox{odd }\alpha:\theta\le\alpha\le\gamma\}|$ and
$R(\theta)=|\{\mbox{even }\alpha:\theta\le\alpha\le\gamma\}|$, and so, by Lemma~\ref{l:binco}(b),
$L(\theta)=R(\theta)$. Now if $\theta=\gamma$, then
$L(\gamma)=|\{\mbox{odd }\alpha:\gamma\le\alpha\le\gamma\}|=1$, and since
$d_\gamma=h(\gamma)$, then $R(\gamma)=1$ as well. The conclusion follows as in Step~1.

\smallbreak
Suppose that $\gamma$ is even.

\smallbreak
3. In order to prove \eq{e:d.7.1}, we first note that if $\gamma=0$, then the double sum
at the right is omitted and the double sum at the left is equal to $h(0)=c_0$. Assume that
$\gamma\ne0$. Then $\LC=\{\mbox{even }\theta:\theta\le\gamma\}$ and
$\RC=\{\gamma\}\cup\{\mbox{even }\theta:%
\exists\,\mbox{odd }\alpha\le\gamma\,\mbox{s.t.}\,\theta\le\alpha\}$ (disjoint union),
and $\LC=\RC$. Let $\theta\in\LC=\RC$. Then
$L(\theta)=|\{\mbox{even }\alpha:\theta\le\alpha\le\gamma\}|$ and, in particular,
$L(\gamma)=1$. If $\theta=\gamma$, then, since $c_\gamma=h(\gamma)$, we have
$R(\gamma)=1$, and if $\theta\ne\gamma$, then
$R(\theta)=|\{\mbox{odd }\alpha:\mbox{$\theta\le\alpha\le\gamma$}\}|$, and so,
by Lemma~\ref{l:binco}(b), $L(\theta)=R(\theta)$.

\smallbreak
4. Finally, we prove \eq{e:d.7.2}. Since the equality is clear for $\gamma=0$
(i.e., `empty' equality), we assume that $|\gamma|>0$. We have
$\LC=\{\mbox{odd }\theta:\theta\le\gamma\}$, $\RC=\{\mbox{odd }\theta:\exists\,%
\mbox{even }\alpha\le\gamma\,\mbox{s.t.}\,\theta\le\alpha\}$ and $\LC=\RC$.
Given $\theta\in\LC=\RC$, we find $\theta\ne\gamma$,
$L(\theta)=|\{\mbox{odd }\alpha:\theta\le\alpha\le\gamma\}|$ and
$R(\theta)=|\{\mbox{even }\alpha:\theta\le\alpha\le\gamma\}|$, and so,
by Lemma~\ref{l:binco}(b), $L(\theta)=R(\theta)$.
\qed\end{pf}

Now we are in a position to prove Theorem~\ref{t:B}.

\begin{prtB}
It suffices to prove only the first inequality: the second one follows from the first
inequality, \eq{e:Vn} and \eq{e:TV}. Setting $u=f(x)$ and $v=f(y)$ and taking into account
\eq{e:fax}, the first inequality in Theorem~\ref{t:B} can be rewritten equivalently as
  \begin{equation} \label{e:eqfB}
d(u,v)\le\sum_{0\ne\alpha\le1}d(u(\alpha),v(\alpha))=
\sum_{j=1}^n\sum_{|\alpha|=j}d(u(\alpha),v(\alpha))
  \end{equation}
(the sum over $|\alpha|=j$ designates the sum over $0\ne\alpha\le1$ s.t.\ $|\alpha|=j$),
where, given $\alpha,\theta\in\mathcal{A}_0$, we set $h(\theta)=f(x+\theta(y-x))$,
  $$u(\alpha)=\sum_{\even\theta\le\alpha}h(\theta),\quad\mbox{and}\quad
    v(\alpha)=\sum_{\odd\theta\le\alpha}h(\theta)\quad\mbox{if $\alpha\ne0$}\quad
    \mbox{and}\quad v(0)\doteq0.$$
In order to establish \eq{e:eqfB}, given integer $0\le j\le n$, we also set
  $$u_j=\sum_{|\alpha|=j}u(\alpha)\qquad\mbox{and}\qquad v_j=\sum_{|\alpha|=j}v(\alpha)$$
and note that
  $$\mbox{$u_0=u(0)=h(0)=u$, $v_0=v(0)\doteq0$, $v=h(1)$, $u_n=u(1)$ and $v_n=v(1)$.}$$
Suppose that we have already verified equalities \eq{e:odm} if $m=n$ is odd and
\eq{e:evm1} if $m=n$ is even. Applying Lemma~\ref{l:metsem}, we get inequality
\eq{e:4d}, where, by virtue of \eq{e:ms.2}, we have:
  $$d(u_j,v_j)=d\biggl(\sum_{|\alpha|=j}u(\alpha),\sum_{|\alpha|=j}v(\alpha)\biggr)\le
    \sum_{|\alpha|=j}d(u(\alpha),v(\alpha)).$$
Now, \eq{e:eqfB} follows if we sum these inequalities over $j=1,\dots,n$ and
take into account \eq{e:4d}.

\smallbreak
It remains to verify equalities \eq{e:odm} and \eq{e:evm1}. For this, we apply
Lemma~\ref{l:funcah} with $\gamma=1$ and note that $m=n=|\gamma|=|1|$. Suppose
that $n=|1|$ is odd. By virtue of \eq{e:d.7.1}, we have:
  \begin{eqnarray}
u+\sum_{i=1}^{(m-1)/2}u_{2i}&=&\sum_{i=0}^{(n-1)/2}u_{2i}=
  \sum_{i=0}^{(n-1)/2}\sum_{|\alpha|=2i}u(\alpha)=
  \sum_{\even\alpha\le1}u(\alpha)\nonumber\\
&=&\sum_{\odd\alpha\le1}u(\alpha)=\sum_{i=1}^{(n+1)/2}\sum_{|\alpha|=2i-1}u(\alpha)=
  \sum_{i=1}^{(m+1)/2}u_{2i-1},\nonumber
  \end{eqnarray}
and by virtue of \eq{e:d.7.2}, we get:
  \begin{eqnarray}
v\!+\!\sum_{i=1}^{(m-1)/2}\!\!v_{2i}&=&h(1)\!+\!\!\sum_{i=0}^{(n-1)/2}\!\!v_{2i}=
  h(1)\!+\!\!\sum_{i=0}^{(n-1)/2}\sum_{|\alpha|=2i}\!\!v(\alpha)=
  h(1)\!+\!\!\sum_{\even\alpha\le1}\!\!v(\alpha)\nonumber\\
&=&\sum_{\odd\alpha\le1}v(\alpha)=\sum_{i=1}^{(n+1)/2}\sum_{|\alpha|=2i-1}v(\alpha)=
  \sum_{i=1}^{(m+1)/2}v_{2i-1},\nonumber
  \end{eqnarray}
which establishes \eq{e:odm}. Now suppose that $n=|1|$ is even. By \eq{e:d.7.1}, we get:
  \begin{eqnarray}
u+\sum_{i=1}^{m/2}u_{2i}&=&\sum_{i=0}^{n/2}u_{2i}=
  \sum_{i=0}^{n/2}\sum_{|\alpha|=2i}u(\alpha)=
  \sum_{\even\alpha\le1}u(\alpha)\nonumber\\
&=&h(1)+\!\sum_{\odd\alpha\le1}\!u(\alpha)=v+\!\sum_{i=1}^{n/2}\sum_{|\alpha|=2i-1}\!\!u(\alpha)=
  v+\!\sum_{i=1}^{m/2}u_{2i-1},\nonumber
  \end{eqnarray}
and by virtue of \eq{e:d.7.2}, we have:
  \begin{eqnarray}
\sum_{i=1}^{m/2}v_{2i}&=&\sum_{i=0}^{n/2}v_{2i}=
  \sum_{i=0}^{n/2}\sum_{|\alpha|=2i}v(\alpha)=
  \sum_{\even\alpha\le1}v(\alpha)\nonumber\\
&=&\sum_{\odd\alpha\le1}\!\!v(\alpha)=\sum_{i=1}^{n/2}\sum_{|\alpha|=2i-1}\!\!v(\alpha)=
  \sum_{i=1}^{m/2}v_{2i-1},\nonumber
  \end{eqnarray}
which establishes \eq{e:evm1} and completes the proof of Theorem~\ref{t:B}.
\qed\end{prtB}

\begin{rem} \label{r:tB}
The left-hand side inequality in Theorem~\ref{t:B} is of interest when $x<y$.
However, if $x\le y$ and $x\not<y$, it can be refined in the following way
(cf.\ \cite[Part~I, Lemma~6]{na05}): given $x,y\in I_a^b$, $x<y$, and $0\ne\gamma\le1$,
we have:
  $$d\bigl(f(x),f(x+\gamma(y-x))\bigr)\le\sum_{0\ne\alpha\le\gamma}%
    \mbox{md}_{|\alpha|}(f_\alpha^x,I_x^y\lfloor\alpha).$$
In fact, by Theorem~\ref{t:B}, we find
  $$d\bigl(f(x),f(x+\gamma(y-x))\bigr)\le\sum_{0\ne\alpha\le1}%
    \mbox{md}_{|\alpha|}(f_\alpha^x,I_x^{x+\gamma(y-x)}\lfloor\alpha),$$
where, by virtue of \eq{e:fax}, the mixed difference at the right is equal to
  \begin{equation} \label{e:d.10}
d\biggl(\sum_{\even\theta\le\alpha}f(x+\theta\gamma(y-x)),
\sum_{\odd\ov\theta\le\alpha}f(x+\ov\theta\gamma(y-x))\biggr).
  \end{equation}
If $\alpha\not\le\gamma$, then $\alpha_i=1$ and $\gamma_i=0$ for some
$i\in\{1,\dots,n\}$, and so, arguing as in Remark~\ref{remA} we find
$x+\theta\gamma(y-x)=x+\ov\theta\gamma(y-x)$ for all even $\theta$ with
$\theta\le\alpha$ implying that \eq{e:d.10} is equal to zero. Now if
$\alpha\le\gamma$, then $\theta\gamma=\theta$ for any $\theta\le\alpha$,
and so, \eq{e:d.10} coincides with the right-hand side of~\eq{e:fax}.
\end{rem}

\section{Proof of Lemma \ref{l:deux}} \label{s:prdeux}

In order to prove Lemma~\ref{l:deux}, we need an auxiliary Lemma~\ref{l:hh},
which plays the same role as Lemma~\ref{l:funcah} above.

\begin{lem} \label{l:hh}
Given a map $h:\mathcal{A}_0\to M$ and a multiindex $\alpha\in\mathcal{A}_0$, we have\/{\rm:}
\newline if $1-\alpha$ is even, then the following two equalities hold\/{\rm:}
  \begin{eqnarray}
\sum_{\even\theta\le\alpha}h(1\!-\!\alpha\!+\!\theta)+\!\!\sum_{\odd\beta\le1-\alpha}\,%
  \sum_{\even\theta\le\alpha+\beta}h(\theta)&=&
  \!\!\sum_{\even\beta\le1-\alpha}\,\sum_{\even\theta\le\alpha+\beta}h(\theta),\label{e:d.11}\\
\sum_{\odd\theta\le\alpha}h(1\!-\!\alpha\!+\!\theta)+\!\!\sum_{\odd\beta\le1-\alpha}\,%
  \sum_{\odd\theta\le\alpha+\beta}h(\theta)&=&
  \!\!\sum_{\even\beta\le1-\alpha}\,\sum_{\odd\theta\le\alpha+\beta}h(\theta),\label{e:d.12}
  \end{eqnarray}
and if $1-\alpha$ is odd, then the following two equalities hold\/{\rm:}
  \begin{eqnarray}
\sum_{1-\alpha\le\even\theta\le1}h(\theta)+\!\!\sum_{\even\beta\le1-\alpha}\,%
  \sum_{\even\theta\le\alpha+\beta}h(\theta)&=&
  \!\!\sum_{\odd\beta\le1-\alpha}\,\sum_{\even\theta\le\alpha+\beta}h(\theta),\label{e:d.13}\\
\sum_{1-\alpha\le\odd\theta\le1}h(\theta)+\!\!\sum_{\even\beta\le1-\alpha}\,%
  \sum_{\odd\theta\le\alpha+\beta}h(\theta)&=&
  \!\!\sum_{\odd\beta\le1-\alpha}\,\sum_{\odd\theta\le\alpha+\beta}h(\theta).\label{e:d.14}
  \end{eqnarray}
\end{lem}

\begin{pf}
As in the proof of Lemma~\ref{l:funcah}, the main idea is to establish equality
\eq{e:LRh}. We divide the proof into four steps.

\smallbreak

Suppose that $1-\alpha$ is even.

\smallbreak
1. Let us prove \eq{e:d.11}. If $\alpha=1$, then $1-\alpha=0$ is even, and equality
\eq{e:d.11} is equivalent to the identity
$\sum_{\even\theta\le1}h(\theta)+0=\sum_{\even\theta\le1}h(\theta)$.
If $\alpha=0$ and if $1-\alpha=1$ is even, then \eq{e:d.11} can be written as
  $$h(1)+\sum_{\odd\beta\le1}\,\sum_{\even\theta\le\beta}h(\theta)=
    \sum_{\even\beta\le1}\,\sum_{\even\theta\le\beta}h(\theta),$$
which was established in \eq{e:d.7.1} for even $\gamma=1$. Thus, in what follows
we assume that $\alpha\ne0,\,1$, i.e., $0<|\alpha|<n$.

We have $\LC=\LC_1\cup\LC_2$, where $\LC_1=\{1-\alpha+\theta':\exists\,\mbox{even }\theta'\le\alpha\}$
(so that $1-\alpha\in\LC_1$) and $\LC_2=\{\mbox{even }\theta:\exists\,\mbox{odd }\beta\le1-\alpha\,\,\,%
\mbox{s.t. }\theta\le\alpha+\beta\}$ (so that $0\in\LC_2$), and
$\RC=\{\mbox{even }\theta:\exists\,\mbox{even }\beta\le1-\alpha\,\,\,\mbox{s.t. }\theta\le\alpha+\beta\}$,
i.e., $\RC=\{\mbox{even }\theta:\theta\le1\}$. We are going to show that $\LC=\RC$.
This equality follows immediately from the definition of $\RC$ and the following
two assertions:
  \begin{eqnarray}
\theta\in\LC_1&\iff&\mbox{$\theta$ is even and $\alpha\lor\theta=1$,}\label{e:d.15}\\[2pt]
\theta\in\LC_2&\iff&\mbox{$\theta$ is even and $\alpha\lor\theta\ne1$,}\label{e:d.16}
  \end{eqnarray}
where $\alpha\lor\theta\equiv\max\{\alpha,\theta\}=\alpha+\theta-\alpha\,\theta$;
in particular, \eq{e:d.15} and \eq{e:d.16} imply that $\LC_1$ and $\LC_2$ are disjoint.
Let us prove \eq{e:d.15}. If $\theta\in\LC_1$, then $\theta=1-\alpha+\theta'$ for
some even $\theta'\le\alpha$ and, since $1-\alpha$ is even and $|\theta|=|1-\alpha|+|\theta'|$,
then $\theta$ is even, $\theta\le(1-\alpha)+\alpha=1$ and
  $$\alpha\lor\theta=\alpha+(1-\alpha+\theta')-\alpha(1-\alpha+\theta')
    =1+\theta'-\alpha\,\theta'=1.$$
Conversely, if $\theta$ is even and $\alpha\lor\theta=1$, then
$\alpha+\theta-\alpha\theta=1$ or $\alpha+\theta=1+\alpha\theta\ge1$.
Setting $\theta'=\alpha+\theta-1$, we find $\theta=1-\alpha+\theta'$, where
$|\theta'|=|\alpha|+|\theta|-n=|\theta|-|1-\alpha|$ is even and $\theta'\le\alpha$,
and so, $\theta\in\LC_1$. Now we establish \eq{e:d.16}.
If $\theta\in\LC_2$, then $\theta$ is even and there exists odd $\beta\le1-\alpha$
s.t.\ $\theta\le\alpha+\beta$, and so, $\alpha\le\alpha+\beta$ and
$\theta\le\alpha+\beta$ imply $\alpha\lor\theta\le\alpha+\beta$. Since
$\beta$ is odd, $1-\alpha$ is even and $\beta\le1-\alpha$, we have
$|\beta|<|1-\alpha|=n-|\alpha|$. It follows that
  $$|\alpha\lor\theta|\le|\alpha+\beta|=|\alpha|+|\beta|<|\alpha|+(n-|\alpha|)=n,$$
and so, $\alpha\lor\theta\ne1$. Conversely, if $\theta$ is even and $\alpha\lor\theta\ne1$,
then there exists $i\in\{1,\dots,n\}$ s.t.\ $\alpha_i=0$ and $\theta_i=0$.
Setting $\beta=(\beta_1,\dots,\beta_n)$ with $\beta_i=0$ and $\beta_j=1-\alpha_j$
if $j\ne i$, we find $\beta\le1-\alpha$, $|\beta|=|1-\alpha|-1$ is odd and
$\theta\le\alpha+\beta$, and so, $\theta\in\LC_2$.

In order to calculate the values $L(\theta)$ and $R(\theta)$ for $\theta\in\LC=\RC$,
we note that, given $0\le\beta\le1-\alpha$, we have:
  \begin{equation} \label{e:at}
\mbox{$\theta\le\alpha+\beta$\quad is equivalent to\quad $(1-\alpha)\theta\le\beta$.}
  \end{equation}
In fact, condition $0\le\beta\le1-\alpha$ is equivalent to condition $\alpha\beta=0$:
  $$0\le\beta\le1-\alpha\iff\beta(1-\alpha)=\beta\iff\beta-\alpha\beta=\beta\iff\alpha\beta=0,$$
and so, if $\theta\le\alpha+\beta$, then $(1-\alpha)\theta\le(1-\alpha)(\alpha+\beta)%
=(1-\alpha)\alpha+\beta-\alpha\beta=\beta$, and if $(1-\alpha)\theta\le\beta$, then
$\theta-\alpha\theta\le\beta$, and so, $\theta\le\alpha\theta+\beta\le\alpha+\beta$.

Given $\theta\in\RC$, by virtue of \eq{e:at}, we find
  $$R(\theta)=|\{\mbox{even }\beta:\mbox{$\beta\!\le\!1\!-\!\alpha$ and $\theta\!\le\!\alpha\!+\!\beta$}\}|
    =|\{\mbox{even }\beta:(1-\alpha)\theta\!\le\!\beta\!\le\!1-\alpha\}|.$$
If $\theta\in\LC_1$, then there exists a unique even $\theta'\le\alpha$ s.t.\
$\theta=1-\alpha+\theta'$, and so, since $\theta\notin\LC_2$, then $L(\theta)=1$.
At the same time,
  $$(1-\alpha)\theta=(1-\alpha)(1-\alpha+\theta')=(1-\alpha)^2+(1-\alpha)\theta'=1-\alpha,$$
and so, by the above, $R(\theta)=1$ as well. Suppose now that $\theta\in\LC_2$. Then,
by \eq{e:d.16}, $1\ne\alpha\lor\theta=\alpha+\theta-\alpha\theta=\alpha+(1-\alpha)\theta$
or $(1-\alpha)\theta\ne1-\alpha$, and so, taking into account \eq{e:at} and
Lemma~\ref{l:binco}(b) we find that
  $$L(\theta)=|\{\mbox{odd }\beta:\mbox{$\beta\!\le\!1\!-\!\alpha$ and $\theta\!\le\!\alpha\!+\!\beta$}\}|
    =|\{\mbox{odd }\beta:(1-\alpha)\theta\!\le\!\beta\!\le\!1-\alpha\}|$$
is equal to $R(\theta)$.

\smallbreak
In the rest of the proof we exhibit only the essential ingredients and differences.

\smallbreak
2. Let us establish \eq{e:d.12}. If $\alpha=1$, we get an identity, and if $\alpha=0$ and
$1=1-\alpha$ is even, we get equality \eq{e:d.7.2} with even $\gamma=1$, and so, we suppose
that $0<|\alpha|<n$. We have $\LC=\LC_1\cup\LC_2$, where
$\LC_1=\{1-\alpha+\theta':\exists\,\mbox{odd }\theta'\le\alpha\}$ and
$\LC_2=\{\mbox{odd }\theta:\exists\,\mbox{odd }\beta\le1-\alpha\,\,\mbox{s.t.\ }\theta\le\alpha+\beta\}$,
and
$\RC=\{\mbox{odd }\theta:\exists\,\mbox{even }\beta\le1-\alpha\,\,\mbox{s.t.\ }\theta\le\alpha+\beta\}$,
which, actually, is $\RC=\{\mbox{odd }\theta:\theta\le1\}$. We need to verify only that
$\LC_1$ and $\LC_2$ are nonempty: the rest of the proof of \eq{e:d.12} (including
\eq{e:d.15} and \eq{e:d.16}) is the same as in Step~1 where `even~$\theta$' is replaced
by `odd~$\theta$'.

Since $\alpha\ne0$, there exists $i\in\{1,\dots,n\}$ s.t.\ $\alpha_i=1$, and so, if we set
$\theta'=(\theta_1',\dots,\theta_n')$ with $\theta_i'=1$ and $\theta_j'=0$ if $j\ne i$, then
$|\theta'|=1$ is odd and $\theta'\le\alpha$. It follows that $1-\alpha+\theta'\in\LC_1$.

Since $\alpha\ne1$, there exists $i\in\{1,\dots,n\}$ s.t.\ $\alpha_i=0$, and so if we
set $\beta=(\beta_1,\dots,\beta_n)$ with $\beta_i=0$ and $\beta_j=1-\alpha_j$ if $j\ne i$,
then $|\beta|=|1-\alpha|-1$ is odd and $\beta\le1-\alpha$. Given $k\in\{1,\dots,n\}$,
$k\ne i$, setting $\theta=(\theta_1,\dots,\theta_n)$ with $\theta_k=1$ and $\theta_j=0$
if $j\ne k$, we find $|\theta|=1$ is odd and $\theta\le\alpha+\beta$, and so, $\theta\in\LC_2$.

\smallbreak
Assume now that $1-\alpha$ is odd. Note that $\alpha\ne1$.

\smallbreak
3. Let us prove \eq{e:d.13}. If $\alpha=0$ and $1=1-\alpha$ is odd, then (since
\mbox{ev\,$\theta$=1} cannot hold in the first sum at the left of \eq{e:d.13})
equality \eq{e:d.13} is equivalent to \eq{e:d.7.1} with odd $\gamma=1$.
Thus, we assume that $|\alpha|>0$.

We have $\LC=\LC_1\cup\LC_2$, where $\LC_1=\{\mbox{even }\theta:1-\alpha\le\theta\le1\}$ and
$\LC_2=\{\mbox{even }\theta:\exists\,\mbox{even }\beta\le1-\alpha\,\,\mbox{s.t.\ }\theta\le\alpha+\beta\}$,
and
$\RC=\{\mbox{even }\theta:\exists\,\mbox{odd }\beta\le1-\alpha\,\,\mbox{s.t.\ }\theta\le\alpha+\beta\}$,
and so, $\RC=\{\mbox{even }\theta:\theta\le1\}$. We have to show that $\LC=\RC$.

First, we show that $\LC_1$ and $\LC_2$ are nonempty. Since $\alpha\ne0$, $\alpha_i=1$
for some $i\in\{1,\dots,n\}$, and so, setting $\theta=(\theta_1,\dots,\theta_n)$ with
$\theta_i=1$ and $\theta_j=1-\alpha_j$ if $j\ne i$, we find that $1-\alpha\le\theta\le1$
and $|\theta|=|1-\alpha|+1$ is even, whence $\theta\in\LC_1$.
Now, since $\alpha\ne1$, $\alpha_i=0$ for some $i\in\{1,\dots,n\}$, and if we set
$\beta=(\beta_1,\dots,\beta_n)$ with $\beta_i=0$ and $\beta_j=1-\alpha_j$ if $j\ne i$,
then we find that $|\beta|=|1-\alpha|-1$ is even, $\theta=0$ is even and
$0\le\alpha+\beta$, and so, $0\in\LC_2$.

Second, we assert that \eq{e:d.15} and \eq{e:d.16} hold; this will imply that $\LC_1$
and $\LC_2$ are disjoint and $\LC=\RC$. In order to prove \eq{e:d.15}, we let
$\theta\in\LC_1$. Then $\theta$ is even and $1-\alpha\le\theta\le1$, and so,
  $$\alpha\lor\theta=\alpha+\theta-\alpha\theta=\alpha+(1-\alpha)\theta=\alpha+(1-\alpha)=1.$$
Conversely, if $\theta$ is even and $\alpha\lor\theta=1$, then $\alpha+\theta-\alpha\theta=1$,
and so, $(1-\alpha)\theta=1-\alpha$ implying $1-\alpha\le\theta$ and $\theta\in\LC_1$.
The proof of \eq{e:d.16} follows the same lines as in Step~1 if `\mbox{odd $\beta$}'
is replaced by `\mbox{even $\beta$}'.

Given $\theta\in\RC$, taking into account \eq{e:at}, we have
$R(\theta)=|\{\mbox{odd }\beta:\mbox{$(1-\alpha)\theta\le\beta\le1-\alpha$}\}|$.
If $\theta\in\LC_1$, then $\theta\notin\LC_2$, and so, $L(\theta)=1$;
in this case $1-\alpha\le\theta$, and so, $(1-\alpha)\theta=1-\alpha$ and $R(\theta)=1$.
Now if $\theta\in\LC_2$, then $\alpha\lor\theta\ne1$, and so,
$(1-\alpha)\theta\ne1-\alpha$ and, by virtue of Lemma~\ref{l:binco}(b), the value
$L(\theta)=|\{\mbox{even }\beta:(1-\alpha)\theta\le\beta\le1-\alpha\}|$ is equal to $R(\theta)$.

\smallbreak
4. Finally, we establish \eq{e:d.14}. If $\alpha=0$ and $1=1-\alpha$ is odd, we get
equality \eq{e:d.7.2} with odd $\gamma=1$. Assume that $|\alpha|>0$. We have
$\LC=\LC_1\cup\LC_2$, where $\LC_1=\{\mbox{odd }\theta:1-\alpha\le\theta\le1\}$
(and so, $1-\alpha\in\LC_1$) and
$\LC_2=\{\mbox{odd }\theta:\exists\,\mbox{even }\beta\le1-\alpha\,\,\mbox{s.t.\ }\theta\le\alpha+\beta\}$,
and
$\RC=\{\mbox{odd }\theta:\exists\,\mbox{odd }\beta\le1-\alpha\,\,\mbox{s.t.\ }\theta\le\alpha+\beta\}$,
and so, $\RC=\{\mbox{odd }\theta:\theta\le1\}$. That $\LC_2$ is nonempty can be seen
as follows. Since $\alpha\ne1$, $\alpha_i=0$ for some $i\in\{1,\dots,n\}$, and so,
if we set $\beta=(\beta_1,\dots,\beta_n)$ with $\beta_i=0$ and $\beta_j=1-\alpha_j$
if $j\ne i$, then $\beta\le1-\alpha$ and $|\beta|=|1-\alpha|-1$ is even. Now, since $\alpha\ne0$,
$\alpha_k=1$ for some $k\ne i$. If we set $\theta=(\theta_1,\dots,\theta_n)$ with
$\theta_k=1$ and $\theta_j=0$ if $j\ne k$, then $|\theta|=1$ is odd and $\theta\le\alpha+\beta$,
and so, $\theta\in\LC_2$. Assertion \eq{e:d.15} with `$\theta$ is even' replaced by
`$\theta$ is odd' is established as in Step~3, while the proof of \eq{e:d.16} follows
the same lines as in Step~1 with `\mbox{odd $\beta$}' replaced by `\mbox{even $\beta$}'.
It follows that $\LC=\RC$. The proof completes with the last paragraph of Step~3.
\qed\end{pf}

\begin{prLdeux}
The inequality (actually, equality) is clear if $\alpha=1$, and so, we assume that $\alpha\ne1$.
The mixed difference at the left-hand side of the inequality is given by \eq{e:fax}, while
given $\alpha\le\beta\le1$, noting that $\alpha\beta=\alpha$ and applying equality \eq{e:2.4*}
we get the following expression for the mixed difference at the right-hand side
(cf.\ \cite[Part~I, expression (3.7)]{na05}):
  $$\mbox{md}_{|\beta|}\bigl(f_\beta^a,I_{a+\alpha(x-a)}^{x+\alpha(y-x)}\lfloor\beta\bigr)
    =d\biggl(\sum_{\even\theta\le\beta}h(\theta),\sum_{\odd\theta\le\beta}h(\theta)\biggr),$$
where $h(\theta)=f\bigl(a+(\alpha\lor\theta)(x-a)+\alpha\theta(y-x)\bigr)$ and
$\alpha\lor\theta=\alpha+\theta-\alpha\theta$.
Changing the summation multiindex $\beta\mapsto\beta-\alpha$ in the sum at the right
of the inequality in Lemma~\ref{l:deux}, we find that it is equivalent to
  $$d(u,v)\le\sum_{0\le\beta\le1-\alpha}d\biggl(\sum_{\even\theta\le\alpha+\beta}h(\theta),
    \sum_{\odd\theta\le\alpha+\beta}h(\theta)\biggr),$$
where
  $$u=\sum_{\even\theta\le\alpha}f(x+\theta(y-x))\qquad\mbox{and}\qquad
    v=\sum_{\odd\theta\le\alpha}f(x+\theta(y-x)).$$
Setting
  $$u(\beta)=\sum_{\even\theta\le\alpha+\beta}h(\theta)\quad\mbox{and}\quad
    v(\beta)=\sum_{\odd\theta\le\alpha+\beta}h(\theta)\quad\mbox{if}\quad0\le\beta\le1-\alpha,$$
the last inequality can be rewritten as
  \begin{equation} \label{e:d.18}
d(u,v)\le\sum_{0\le\beta\le1-\alpha}d(u(\beta),v(\beta))=\sum_{j=0}^{|1-\alpha|}%
\sum_{|\beta|=j}d(u(\beta),v(\beta)).
  \end{equation}

In order to establish \eq{e:d.18}, we will apply Lemma~\ref{l:metsem} with
$m=|1-\alpha|+1=n-|\alpha|+1$ and
  $$u_j=\sum_{|\beta|=j-1}u(\beta)\quad\mbox{and}\quad v_j=\sum_{|\beta|=j-1}v(\beta)
    \quad\mbox{if}\quad 1\le j\le m.$$
Suppose that we have already verified equalities \eq{e:odm} and \eq{e:evm2}. Then
by Lemma~\ref{l:metsem}, we get inequality \eq{e:4d}, where, by virtue of \eq{e:ms.2},
  $$d(u_j,v_j)=d\biggl(\sum_{|\beta|=j-1}u(\beta),\sum_{|\beta|=j-1}v(\beta)\biggr)\le
    \sum_{|\beta|=j-1}d(u(\beta),v(\beta)),\quad 1\le j\le m.$$
Summing over $j=1,\dots,m$ and taking into account \eq{e:4d}, we arrive at \eq{e:d.18}:
  $$d(u,v)\le\sum_{j=1}^md(u_j,v_j)\le\sum_{j=1}^{|1-\alpha|+1}\sum_{|\beta|=j-1}d(u(\beta),v(\beta)).$$

\smallbreak
Assume that $1-\alpha$ is even; then $m$ is odd. Let us verify the first equality
in \eq{e:odm}. For this, we apply equality \eq{e:d.11} and calculate the first sum
at the left-hand side of \eq{e:d.11}. Given even $\theta\le\alpha$, we have
$1-\alpha+\theta\in\LC_1$ (cf.\ Step~1 in the proof of Lemma~\ref{l:hh}), and so,
by \eq{e:d.15}, $\alpha\lor(1-\alpha+\theta)=1$ and $\alpha(1-\alpha+\theta)=\theta$,
so that the definition of $h(1-\alpha+\theta)$ implies
  $$\sum_{\even\theta\le\alpha}h(1-\alpha+\theta)=\sum_{\even\theta\le\alpha}f(x+\theta(y-x))=u.$$
Applying equality \eq{e:d.11}, we get:
  \begin{eqnarray}
u+\sum_{i=1}^{(m-1)/2}u_{2i}&=&u+\sum_{i=1}^{|1-\alpha|/2}\sum_{|\beta|=2i-1}u(\beta)
  =u+\sum_{\odd\beta\le1-\alpha}u(\beta)\nonumber\\
&=&\sum_{\even\beta\le1-\alpha}u(\beta)=\sum_{i=0}^{|1-\alpha|/2}\sum_{|\beta|=2i}u(\beta)
  =\sum_{i=0}^{|1-\alpha|/2}u_{2i+1}\nonumber\\
&=&\sum_{i=1}^{(|1-\alpha|+2)/2}u_{2i-1}=\sum_{i=1}^{(m+1)/2}u_{2i-1},\nonumber
  \end{eqnarray}
and the first equality in \eq{e:odm} follows. In a similar manner we find that the
first sum at the left-hand side of \eq{e:d.12} is equal to $v$, and, by virtue of
\eq{e:d.12}, the calculations above show that the second equality in \eq{e:odm}
holds as well.

\smallbreak
Now suppose that $1-\alpha$ is odd, and so, $m$ (defined above) is even.
In order to verify the first equality in \eq{e:evm2}, we calculate the first
sum at the left-hand side of \eq{e:d.13}. Given even $\theta$ with $1-\alpha\le\theta\le1$,
we have (cf.\ Step~3 in the proof of Lemma~\ref{l:hh}) $\theta\in\LC_1$ and
$\alpha\lor\theta=1$. Moreover (cf.\ \cite[Part~I, assertion (3.9)]{na05}),
there exists a unique $\theta'\in\mathcal{A}_0$ s.t.\ $\theta'\le\alpha$ and
$\theta=1-\alpha+\theta'$ (define $\theta'$ by $\theta'=\alpha+\theta-1$).
Since $|\theta'|=|\alpha|+|\theta|-n=|\theta|-|1-\alpha|$ and $1-\alpha$ is odd,
then $\theta'$ is odd, and $\alpha\theta=\alpha(1-\alpha+\theta')=\theta'$.
It follows that $h(\theta)=f(x+\theta'(y-x))$. Changing the summation multiindex
$\theta\mapsto\theta'$ in the first sum at the left of \eq{e:d.13}, we get:
  $$\sum_{1-\alpha\le\even\theta\le1}h(\theta)=\sum_{\odd\theta'\le\alpha}%
    f(x+\theta'(y-x))=v.$$
Applying equality \eq{e:d.13}, we find
  \begin{eqnarray}
\sum_{i=1}^{m/2}u_{2i}&=&\sum_{i=1}^{(|1-\alpha|+1)/2}\sum_{|\beta|=2i-1}u(\beta)
  =\sum_{\odd\beta\le1-\alpha}u(\beta)\nonumber\\
&=&v+\sum_{\even\beta\le1-\alpha}u(\beta)=v+\sum_{i=0}^{(|1-\alpha|-1)/2}\sum_{|\beta|=2i}u(\beta)
  =v+\sum_{i=0}^{(|1-\alpha|-1)/2}u_{2i+1}\nonumber\\
&=&v+\sum_{i=1}^{(|1-\alpha|+1)/2}u_{2i-1}=v+\sum_{i=1}^{m/2}u_{2i-1},\nonumber
  \end{eqnarray}
which proves the first equality in \eq{e:evm2}. Similarly, the first sum at the left-hand
side of \eq{e:d.14} is equal to $u$, and, by virtue of \eq{e:d.14}, the calculations
above prove the second equality in \eq{e:evm2}.

This completes the proof of Lemma~\ref{l:deux}.
\qed\end{prLdeux}

\section{Proof of Lemma \ref{l:troix}} \label{s:prtroix}

Note that if $\mathcal{P}=\{x[\sigma]\}_{\sigma=0}^\kappa$ is a net partition of $I_a^b$, then
  \begin{equation} \label{e:E.1}
I_a^b=\bigcup_{1\le\sigma\le\kappa}I_{x[\sigma-1]}^{x[\sigma]}=
\bigcup_{1\le\sigma\le\kappa}\prod_{i=1}^n[x_i(\sigma_i-1),x_i(\sigma_i)]=
\prod_{i=1}^n\biggl(\bigcup_{l=1}^{\kappa_i}I_{x_i(l-1)}^{x_i(l)}\biggr)
  \end{equation}
is a union of non-overlapping non-degenerated rectangles $I_{x[\sigma-1]}^{x[\sigma]}$ with
the sides parallel to the coordinate axes. In this section it will be convenient and brief to
term the union as in \eq{e:E.1} also a {\em partition} of $I_a^b$ (by non-overlapping
non-degenerated rectangles).

\smallbreak
If $\mathcal{P}=\{x[\sigma]\}_{\sigma=0}^\kappa$ and
$\mathcal{P}'=\{x'[\sigma']\}_{\sigma'=0}^{\kappa'}$ are two net partitions of $I_a^b$,
we say that $\mathcal{P}'$ is a {\em refinement\/} of $\mathcal{P}$ if
$\mathcal{P}\subset\mathcal{P}'$. Also, for the sake of convenience we define the
{\em $n$-th prevariation\/} of $f:I_a^b\to M$, corresponding to $\mathcal{P}$, by
  $$\mbox{\rm v}_n(f;\mathcal{P})=\sum_{1\le\sigma\le\kappa}\mbox{md}_n(f,I_{x[\sigma-1]}^{x[\sigma]}).$$
It follows that the Vitali-type $n$-th variation of $f$ is given by
$V_n(f,I_a^b)=\sup_{\mathcal{P}}\mbox{\rm v}_n(f;\mathcal{P})$, where the supremum
is taken over all net partitions $\mathcal{P}$ of $I_a^b$.

\smallbreak
The basic ingredient in the proof of Lemma~\ref{l:troix} is the following

\begin{lem} \label{l:Eone}
Given $f:I_a^b\to M$, if $\mathcal{P}$ and $\mathcal{P}'$ are two net partitions
of $I_a^b$ s.t.\ $\mathcal{P}\subset\mathcal{P}'$, then
$\mbox{\rm v}_n(f;\mathcal{P})\le\mbox{\rm v}_n(f;\mathcal{P}')$.
\end{lem}

In order to prove this lemma we need three more Lemmas~\ref{l:Etwo}--\ref{l:Efour}.
In what follows we fix a map $f:I_a^b\to M$.

\begin{lem} \label{l:Etwo}
Given $x,y\in I_a^b$ with $x<y$ and $x'\in I_a^b$, we have the following {\sl partition
of $I_x^y$, induced by the point $x'$\/{\rm:}}
  \begin{equation} \label{e:E.2}
I_x^y=\bigcup_{1-\xi\le\alpha\le1}I_{x+\alpha\xi(x'-x)}^{x'+\alpha(y-x')},
  \end{equation}
where the multiindex $\xi\equiv\xi(x,x',y)=(\xi_1,\dots,\xi_n)$ is given by
  \begin{equation} \label{e:E.3}
\xi_i\equiv\xi_i(x,x',y)=\left\{
  \begin{array}{rcl}
  1 &\mbox{if}& x_i<x_i'<y_i,\\[2pt]
  0 &\mbox{if}& \mbox{$x_i'\le x_i$ or $x_i'\ge y_i$},
  \end{array}\right.
  \quad\,\,i\in\{1,\dots,n\},
  \end{equation}
and
  \begin{equation} \label{e:E.4}
\mbox{\rm md}_n(f,I_x^y)\le\sum_{1-\xi\le\alpha\le1}%
\mbox{\rm md}_n\bigl(f,I_{x+\alpha\xi(x'-x)}^{x'+\alpha(y-x')}\bigr).
  \end{equation}
\end{lem}

Before we prove Lemma~\ref{l:Etwo}, let us establish two of its particular variants
as Lemmas \ref{l:Ethree} and \ref{l:Efour} (note that in Lemma~\ref{l:Ethree} the
rectangles in the union may degenerate).

\begin{lem} \label{l:Ethree}
If $x,y\in I_a^b$ with $x<y$ and $x'\in I_x^y$, then we have the following union of
non-overlapping\/ {\rm(}possibly, degenerated\/{\rm)} rectangles
  \begin{equation} \label{e:E.4,5}
I_x^y=\bigcup_{0\le\alpha\le1}I_{x+\alpha(x'-x)}^{x'+\alpha(y-x')},
  \end{equation}
and the following inequality holds\/{\rm:}
  \begin{equation} \label{e:E.5}
\mbox{\rm md}_n(f,I_x^y)\le\sum_{0\le\alpha\le1}%
\mbox{\rm md}_n\bigl(f,I_{x+\alpha(x'-x)}^{x'+\alpha(y-x')}\bigr).
  \end{equation}
\end{lem}

\begin{pf}
Since $x_i\le x_i'\le y_i$ for all $i\in\{1,\dots,n\}$, we have:
  $$I_{x_i}^{y_i}=[x_i,y_i]=[x_i,x_i']\cup[x_i',y_i]=I_{x_i}^{x_i'}\cup I_{x_i'}^{y_i}
    =\bigcup_{\alpha_i=0}^1I_{x_i+\alpha_i(x_i'-x_i)}^{x_i'+\alpha_i(y_i-x_i')},$$
and so (cf.\ equation (2.5) in \cite[Part~II]{na05}),
  $$I_x^y=\prod_{i=1}^nI_{x_i}^{y_i}=\prod_{i=1}^n\biggl(%
    \bigcup_{\alpha_i=0}^1I_{x_i+\alpha_i(x_i'-x_i)}^{x_i'+\alpha_i(y_i-x_i')}\biggr)
    =\bigcup_{0\le\alpha\le1}I_{x+\alpha(x'-x)}^{x'+\alpha(y-x')}.$$

The mixed difference at the left-hand side of \eq{e:E.5} is given by \eq{e:mdn},
and again by virtue of \eq{e:mdn}, the mixed difference at the right-hand side
of \eq{e:E.5} is equal to
  $$\mbox{md}_n(f,I_{x+\alpha(x'-x)}^{x'+\alpha(y-x')})=d\biggl(
    \sum_{\even\beta\le1}h(\alpha,\beta),\sum_{\odd\beta\le1}h(\alpha,\beta)\biggr),$$
where $h(\alpha,\beta)=f(x+(\alpha\lor\beta)(x'-x)+\alpha\beta(y-x'))$ and
$\alpha\lor\beta=\alpha+\beta-\alpha\beta$. Noting that if $\alpha=\beta$, then
$\alpha\lor\beta=\beta$ and $\alpha\beta=\beta$, we find
  \begin{eqnarray}
\sum_{0\le\alpha\le1}\,\sum_{\even\beta\le1}h(\alpha,\beta)&=&
  \sum_{\even\beta\le1}\sum_{%
  \begin{array}{c}\\[-20pt]\scriptstyle 0\le\alpha\le1,\\[-6pt]\scriptstyle\alpha\ne\beta\end{array}}
  h(\alpha,\beta)+\sum_{\even\beta\le1}h(\beta,\beta)\nonumber\\
&=&\sum_{\even\beta\le1}\sum_{%
  \begin{array}{c}\\[-20pt]\scriptstyle 0\le\alpha\le1,\\[-6pt]\scriptstyle\alpha\ne\beta\end{array}}
  h(\alpha,\beta)+\sum_{\even\beta\le1}f(x+\beta(y-x))\nonumber\\
&\equiv&U+u.\nonumber
  \end{eqnarray}
Let us show that the double sum $U$ can be represented as
  $$U=\sum_{0\ne\gamma\le1}\sum_{%
  \begin{array}{c}\\[-20pt]\scriptstyle 0\le\delta\le\gamma,\\[-6pt]\scriptstyle\delta\ne\gamma\end{array}}
  c_{\gamma\delta}f\bigl(x+\gamma(x'-x)+\delta(y-x')\bigr)$$
with certian integer factors $c_{\gamma\delta}$ to be evaluated below. In fact,
given $0\ne\gamma\le1$ and $0\le\delta\le\gamma$ with $\delta\ne\gamma$, there exist
even $\beta\le1$ and $0\le\alpha\le1$, $\alpha\ne\beta$, s.t.\ $\alpha\lor\beta=\gamma$
and $\alpha\beta=\delta$. In order to see this, if $\gamma$ is even or $\delta$ is even,
we may set $\beta=\gamma$ and $\alpha=\delta$, or $\beta=\delta$ and $\alpha=\gamma$,
respectively. Now, if $\gamma$ and $\delta$ are odd, then since $\delta\ne\gamma$,
we can find $i\in\{1,\dots,n\}$ s.t.\ $\delta_i=0$ and $\gamma_i=1$, and so,
if we set $\beta=(\delta_1,\dots,\delta_{i-1},1,\delta_{i+1},\dots,\delta_n)$,
then $\delta\le\beta\le\gamma$, $\delta\ne\beta\ne\gamma$ and $|\beta|=|\delta|+1$
is even, and it remains to put $\alpha=\gamma+\delta-\beta$.

Given $\gamma$ and $\delta$ as above, let us evaluate $c_{\gamma\delta}$. Since
$\delta=\alpha\beta\le\beta\le\alpha\lor\beta=\gamma$ and, given even $\beta$,
the multiindex $0\le\alpha\le1$, $\alpha\ne\beta$, s.t.\ $\alpha\lor\beta=\gamma$
and $\alpha\beta=\delta$, is determined uniquely by $\alpha=\gamma+\delta-\beta$,
we have $c_{\gamma\delta}=|\{\mbox{even }\beta:\delta\le\beta\le\gamma\}|$.

In a similar manner, we find
  $$\sum_{0\le\alpha\le1}\,\sum_{\odd\beta\le1}h(\alpha,\beta)=
    \sum_{\odd\beta\le1}\sum_{%
  \begin{array}{c}\\[-20pt]\scriptstyle 0\le\alpha\le1,\\[-6pt]\scriptstyle\alpha\ne\beta\end{array}}
  h(\alpha,\beta)+\sum_{\odd\beta\le1}f(x+\beta(y-x))\equiv V+v,$$
where
  $$V=\sum_{0\ne\gamma\le1}\sum_{%
  \begin{array}{c}\\[-20pt]\scriptstyle 0\le\delta\le\gamma,\\[-6pt]\scriptstyle\delta\ne\gamma\end{array}}
  d_{\gamma\delta}f\bigl(x+\gamma(x'-x)+\delta(y-x')\bigr)$$
with $d_{\gamma\delta}=|\{\mbox{odd }\beta:\delta\le\beta\le\gamma\}|$.
By Lemma~\ref{l:binco}(b), $c_{\gamma\delta}=d_{\gamma\delta}$, and so, $U=V$.
Applying the translation invariance of $d$ and inequality \eq{e:ms.2},
we obtain inequality \eq{e:E.5}:
  \begin{eqnarray}
d(u,v)&=&d(U+u,V+v)=d\biggl(\sum_{0\le\alpha\le1}\sum_{\even\beta\le1}h(\alpha,\beta),
  \sum_{0\le\alpha\le1}\sum_{\odd\beta\le1}h(\alpha,\beta)\biggr)\nonumber\\
&\le&\sum_{0\le\alpha\le1}d\biggl(\sum_{\even\beta\le1}h(\alpha,\beta),
  \sum_{\odd\beta\le1}h(\alpha,\beta)\biggr).\qquad\qquad\square\nonumber
  \end{eqnarray}
\end{pf}

\begin{rem} \label{r:2t}
If $x<x'<y$ in Lemma~\ref{l:Ethree}, then all rectangles at the right-hand side of
\eq{e:E.4,5} are non-degenerated, i.e., $x+\alpha(x'-x)<x'+\alpha(y-x')$ for all
$0\le\alpha\le1$. Moreover, the point $x'$ gives rise to a net partition
$\{x[\sigma]\}_{\sigma=0}^\kappa$ of $I_x^y$ with $x[\sigma]=(x_1(\sigma_1),\dots,x_n(\sigma_n))$
and $0\le\sigma\le\kappa$ as follows: we put $\kappa=2=1+1\in\mathbb{N}^n$ and,
given $i\in\{1,\dots,n\}$, we set $x_i(0)=x_i$, $x_i(1)=x_i'$ and $x_i(2)=y_i$.
We note that if $0\le\sigma\le1$, then $x[\sigma]=x+\sigma(x'-x)$, and if
$1\le\sigma\le2$, then $x[\sigma]=x'+(\sigma-1)(y-x')$. It follows that
  $$I_x^y=\bigcup_{0\le\alpha\le1}I_{x+\alpha(x'-x)}^{x'+\alpha(y-x')}
    =\bigcup_{1\le\sigma\le2}I_{x+(\sigma-1)(x'-x)}^{x'+(\sigma-1)(y-x')}
    =\bigcup_{1\le\sigma\le\kappa}I_{x[\sigma-1]}^{x[\sigma]}.$$
However, in the general case $x\le x'\le y$ we may also have $x\not<x'$ or
$x'\not<y$, and so, there exists $i\in\{1,\dots,n\}$ s.t.\ $x_i=x_i'$ or $x_i'=y_i$.
Thus, since some coordinates of $x'$ may be equal to the corresponding coordinates
of $x$ and/or $y$, certain rectangles at the right-hand side of \eq{e:E.4,5} may
degenerate into lower-dimensional rectangles, and so, by Remark~\ref{remA}, the
mixed difference $\mbox{md}_n$ over these rectangles is equal to zero.
In order to exclude these degenerated rectangles from the consideration, we establish
the following lemma.
\end{rem}

\begin{lem} \label{l:Efour}
Given $x,y\in I_a^b$ with $x<y$ and $x'\in I_x^y$, we have the following {\sl partition}
of $I_x^y$, induced by $x'$:
  \begin{equation} \label{e:E.6}
I_x^y=\bigcup_{\lambda\le\alpha\le\mu}I_{x+\alpha(x'-x)}^{x'+\alpha(y-x')},
  \end{equation}
where the multiindices $\lambda\equiv\lambda(x,x')=(\lambda_1,\dots,\lambda_n)$ and
$\mu\equiv\mu(x',y)=(\mu_1,\dots,\mu_n)$ are defined for $i\in\{1,\dots,n\}$ by
  $$\lambda_i\equiv\lambda_i(x,x')=\left\{
    \begin{array}{rcl}
    1 &\mbox{if}& x_i=x_i',\\[2pt]
    0 &\mbox{if}& x_i<x_i',
    \end{array}\right.
    \,\,\,\mbox{and}\,\,\,
  \mu_i\equiv\mu_i(x',y)=\left\{
    \begin{array}{rcl}
    0 &\mbox{if}& x_i'=y_i,\\[2pt]
    1 &\mbox{if}& x_i'<y_i,
    \end{array}\right.$$
and the following inequality holds\/{\rm:}
  \begin{equation} \label{e:E.7}
\mbox{\rm md}_n(f,I_x^y)\le\sum_{\lambda\le\alpha\le\mu}\mbox{\rm md}_n%
\bigl(f,I_{x+\alpha(x'-x)}^{x'+\alpha(y-x')}\bigr).
  \end{equation}
\end{lem}

\begin{pf}
First, we note that, since $x_i<y_i$ for all $i\in\{1,\dots,n\}$, then $\lambda\le\mu$.
In particular, if $x<x'<y$, then $\lambda=0$ and $\mu=1$, and we get \eq{e:E.4,5}
as a consequence of \eq{e:E.6}; cf.\ Remark~\ref{r:2t}.

In order to prove \eq{e:E.6}, given $i\in\{1,\dots,n\}$, consider the following
possibilities: (i) $x_i'=x_i$ and $x_i'<y_i$; (ii) $x_i<x_i'$ and $x_i'=y_i$;
and (iii) $x_i<x_i'$ and $x_i'<y_i$. We have, respectively:

(i) $\lambda_i=1$ and $\mu_i=1$, and so, if $\lambda_i\le\alpha_i\le\mu_i$, then $\alpha_i=1$ and
  $$I_{x_i}^{y_i}=I_{x_i'}^{y_i}=\bigcup_{\alpha_i=1}
    I_{x_i+\alpha_i(x_i'-x_i)}^{x_i'+\alpha_i(y_i-x_i')};$$

(ii) $\lambda_i=0$ and $\mu_i=0$, and so, if $\lambda_i\le\alpha_i\le\mu_i$, then $\alpha_i=0$ and
  $$I_{x_i}^{y_i}=I_{x_i}^{x_i'}=\bigcup_{\alpha_i=0}
    I_{x_i+\alpha_i(x_i'-x_i)}^{x_i'+\alpha_i(y_i-x_i')};$$

(iii) $\lambda_i=0$ and $\mu_i=1$, and so, if $\lambda_i\le\alpha_i\le\mu_i$,
then $\alpha_i\in\{0,1\}$ and
  $$I_{x_i}^{y_i}=I_{x_i}^{x_i'}\cup I_{x_i'}^{y_i}=\bigcup_{\alpha_i=0}^1
    I_{x_i+\alpha_i(x_i'-x_i)}^{x_i'+\alpha_i(y_i-x_i')}.$$
Moreover, in all the cases (i)--(iii) the left endpoint $x_i+\alpha_i(x_i'-x_i)$ is
less than the right endpoint $x_i'+\alpha_i(y_i-x_i')$, and so, all the closed
intervals above are non-degenerated. It follows that
  $$I_x^y=\prod_{i=1}^nI_{x_i}^{y_i}=\prod_{i=1}^n\biggl(\bigcup_{\lambda_i\le\alpha_i\le\mu_i}
    I_{x_i+\alpha_i(x_i'-x_i)}^{x_i'+\alpha_i(y_i-x_i')}\biggr)=
    \bigcup_{\lambda\le\alpha\le\mu}I_{x+\alpha(x'-x)}^{x'+\alpha(y-x')}.$$

The point $x'$ gives rise to a net partition $\{x[\sigma]\}_{\sigma=0}^\kappa$ of $I_x^y$
as follows: we put $\kappa=\mu-\lambda+1$ and, given $i\in\{1,\dots,n\}$, we set
$x_i(0)=x_i$ and $x_i(1)=y_i$ if $\kappa_i=1$, and $x_i(0)=x_i$, $x_i(1)=x_i'$ and
$x_i(2)=y_i$ if $\kappa_i=2$. We note that if $0\le\sigma\le\mu-\lambda$, then
$x[\sigma]=x+(\sigma+\lambda)(x'-x)$, and if $1\le\sigma\le\kappa=\mu-\lambda+1$, then
$x[\sigma]=x'+(\sigma-1+\lambda)(y-x')$. Also, note that $x+\lambda(x'-x)=x$ and
$x'+\mu(y-x')=y$. It follows that
  $$I_x^y=\bigcup_{\lambda\le\alpha\le\mu}I_{x+\alpha(x'-x)}^{x'+\alpha(y-x')}=
    \bigcup_{1\le\sigma\le\mu-\lambda+1}I_{x+(\sigma-1+\lambda)(x'-x)}^{x'+(\sigma-1+\lambda)(y-x')}
    =\bigcup_{1\le\sigma\le\kappa}I_{x[\sigma-1]}^{x[\sigma]}.$$

Now, we turn to the proof of \eq{e:E.7}. By Lemma~\ref{l:Ethree}, inequality \eq{e:E.5} holds.
Clearly, if $\lambda=0$ and $\mu=1$ (i.e., $x<x'<y$), then \eq{e:E.5} implies \eq{e:E.7}.
Assume that $\lambda\ne0$ (i.e., $x\not<x'$) and suppose that $0\le\alpha\le1$ is s.t.\
$\lambda\not\le\alpha$. Then there exists $i\in\{1,\dots,n\}$ s.t.\ $\lambda_i=1$
and $\alpha_i=0$, and so, $x_i=x_i'$, which implies
$x_i+\alpha_i(x_i'-x_i)=x_i=x_i'=x_i'+\alpha_i(y_i-x_i')$. It follows from Remark~\ref{remA}
that $\mbox{md}_n(f,I_{x+\alpha(x'-x)}^{x'+\alpha(y-x')})=0$. Similarly, if we assume that
$\mu\ne1$ (i.e., $x'\not<y$) and suppose that $0\le\alpha\le1$ is s.t.\
$\alpha\not\le\mu$, then there exists $i\in\{1,\dots,n\}$ s.t.\ $\alpha_i=1$ and
$\mu_i=0$, and so, $x_i'=y_i$. Noting that
$x_i+\alpha_i(x_i'-x_i)=x_i'=y_i=x_i'+\alpha_i(y_i-x_i')$, we find
$\mbox{md}_n(f,I_{x+\alpha(x'-x)}^{x'+\alpha(y-x')})=0$.
In this way inequality \eq{e:E.7} follows.
\qed\end{pf}

\begin{prEtwo}
Suppose that $x,y\in I_a^b$, $x<y$ and $x'\in I_a^b$. We set $x''=x+\xi(x'-x)$, where
$\xi$ is defined in \eq{e:E.3} (the point $x''$ will play the role of $x'$ from \eq{e:E.6}).
We have $x\le x''<y$; in fact, given $i\in\{1,\dots,n\}$, we find: if $\xi_i=1$, then
$x_i<x_i'<y_i$ and $x_i''=x_i'$ implying $x_i<x_i''<y_i$, and if $\xi_i=0$, then
$x_i'\le x_i$ or $x_i'\ge y_i$, and $x_i''=x_i$ implying $x_i=x_i''<y_i$. Applying
\eq{e:E.6} with $x'$ replaced by $x''$, we get the following partition of $I_x^y$
induced by $x''$ and, hence, by~$x'$:
  \begin{equation} \label{e:E.8}
I_x^y=\bigcup_{\lambda''\le\alpha\le\mu''}I_{x+\alpha(x''-x)}^{x''+\alpha(y-x'')},
  \end{equation}
where $\lambda''=\lambda(x,x'')$ and $\mu''=\mu(x'',y)$ are defined in Lemma~\ref{l:Efour},
i.e., given $i\in\{1,\dots,n\}$, we have:
  $$\lambda_i''=\left\{
    \begin{array}{rcl}
    1 &\mbox{if}& x_i=x_i'',\\[2pt]
    0 &\mbox{if}& x_i<x_i'',
    \end{array}\right.
    \quad\,\mbox{and}\quad\,
  \mu_i''=\left\{
    \begin{array}{rcl}
    0 &\mbox{if}& x_i''=y_i,\\[2pt]
    1 &\mbox{if}& x_i''<y_i.
    \end{array}\right.$$

We assert that $\lambda''=1-\xi$ and $\mu''=1$. In fact, since $x''<y$, then $\mu''=1$.
In order to see that $\lambda''=1-\xi$, let $i\in\{1,\dots,n\}$.
If $x_i<x_i'<y_i$, then $\xi_i=1$, and so, $x_i''=x_i+\xi_i(x_i'-x_i)=x_i'$, which
implies $x_i<x_i''$ and $\lambda_i''=0=1-\xi_i$. Now if $x_i'\le x_i$ or $x_i'\ge y_i$,
then $\xi_i=0$, and so, $x_i''=x_i$, which gives $\lambda_i''=1=1-\xi_i$.

Now, let us calculate the lower and upper indices in \eq{e:E.8}. We have:
$x\!+\!\alpha(x''\!-\!x)=x+\alpha\xi(x'-x)$ and
  $$x''+\alpha(y-x'')=x+(1-\alpha)\xi(x'-x)+\alpha(y-x).$$
Noting that the union in \eq{e:E.8} is taken over $\alpha\le1$ s.t.\ $1-\xi\le\alpha$,
we get $1-\alpha\le\xi$, and so, $(1-\alpha)\xi=1-\alpha$ implying
  $$x''+\alpha(y-x'')=x+(1-\alpha)(x'-x)+\alpha(y-x)=x'+\alpha(y-x').$$
These calculations and observations above prove equality \eq{e:E.2}.

Let us show that partition \eq{e:E.2} is actually induced by $x'$.
Since $x'\in I_a^b$, by Lemma~\ref{l:Efour}, the point $x'$ induces a partition
of $I_a^b$ of the form~\eq{e:E.6}:
  $$I_a^b=\bigcup_{\lambda'\le\beta\le\mu'}I_{a+\beta(x'-a)}^{x'+\beta(b-x')},$$
where the multiindices $\lambda'\!=\!\lambda(a,x')$ and $\mu'\!=\!\mu(x',b)$ are defined
in Lemma~\ref{l:Efour}, i.e., given $i\in\{1,\dots,n\}$, we have:
  $$\lambda_i'=\left\{
    \begin{array}{rcl}
    1 &\mbox{if}& a_i=x_i',\\[2pt]
    0 &\mbox{if}& a_i<x_i',
    \end{array}\right.
    \quad\,\mbox{and}\quad\,
  \mu_i'=\left\{
    \begin{array}{rcl}
    0 &\mbox{if}& x_i'=b_i,\\[2pt]
    1 &\mbox{if}& x_i'<b_i.
    \end{array}\right.$$

We assert that for each $\alpha$ with $1-\xi\le\alpha\le1$ there exists a unique
$\beta\equiv\beta(\alpha)$ with $\lambda'\le\beta\le\mu'$ s.t.\
  \begin{equation} \label{e:E.9}
I_{x+\alpha\xi(x'-x)}^{x'+\alpha(y-x')}=I_x^y\cap I_{a+\beta(x'-a)}^{x'+\beta(b-x')}.
  \end{equation}
In order to prove \eq{e:E.9}, we define $\beta=\beta(\alpha)=(\beta_1,\dots,\beta_n)$ by
  $$\beta_i\equiv\beta_i(\alpha)=\left\{
    \begin{array}{ccl}
    \alpha_i &\mbox{if}& x_i'<y_i,\\[2pt]
    0 &\mbox{if}& x_i'\ge y_i,
    \end{array}\right.
    \qquad i\in\{1,\dots,n\},$$
and establish equality \eq{e:E.9} componentwise. Given $i\in\{1,\dots,n\}$,
we consider the following two cases: (a) $x_i'<y_i$, and (b) $x_i'\ge y_i$.

In case (a) we have $\beta_i=\alpha_i$. First, assume that $x_i<x_i'$, and so,
$\xi_i=1$. It follows that if $1-\xi_i\le\alpha_i\le1$, then $\alpha_i=0$ or
$\alpha_i=1$. If $\alpha_i=0$, then we find (for $\beta_i=\alpha_i=0$)
  $$I_{x_i}^{x_i'}=[x_i,x_i']=[x_i,y_i]\cap[a_i,x_i']=
    I_{x_i}^{y_i}\cap I_{a_i+\beta_i(x_i'-a_i)}^{x_i'+\beta_i(b_i-x_i')},$$
and if $\alpha_i=1$, then we find (for $\beta_i=\alpha_i=1$)
  $$I_{x_i'}^{y_i}=[x_i',y_i]=[x_i,y_i]\cap[x_i',b_i]=
    I_{x_i}^{y_i}\cap I_{a_i+\beta_i(x_i'-a_i)}^{x_i'+\beta_i(b_i-x_i')}.$$
Now, assume that $x_i'\le x_i$, and so, $\xi_i=0$ and $x_i'\le x_i<y_i\le b_i$.
It follows that if $1-\xi_i\le\alpha_i\le1$, then $\beta_i=\alpha_i=1$ implying
  $$I_{x_i}^{y_i}=[x_i,y_i]=[x_i,y_i]\cap[x_i',b_i]=
    I_{x_i}^{y_i}\cap I_{a_i+\beta_i(x_i'-a_i)}^{x_i'+\beta_i(b_i-x_i')}.$$

In case (b) we have $\xi_i=0$, $\beta_i=0$ and $a_i\le x_i<y_i\le x_i'$, and so,
if $1-\xi_i\le\alpha_i\le1$, then $\alpha_i=1$ and
  $$I_{x_i}^{y_i}=[x_i,y_i]=[x_i,y_i]\cap[a_i,x_i']=
    I_{x_i}^{y_i}\cap I_{a_i+\beta_i(x_i'-a_i)}^{x_i'+\beta_i(b_i-x_i')}.$$

Let us show that $\lambda'\le\beta\le\mu'$. Let $i\in\{1,\dots,n\}$. If $a_i=x_i'$,
then $\lambda_i'=1=\mu_i'$ and, since $x_i'<y_i$, then $\beta_i=\alpha_i$.
By \eq{e:E.3}, $\xi_i=0$, and so, since $1-\xi_i\le\alpha_i\le1$, then $\alpha_i=1$,
which implies $\lambda_i'=\beta_i=\mu_i'$. Now, if $x_i'=b_i$, then
$\lambda_i'=0=\mu_i'$ and, since $x_i'\ge y_i$, then $\beta_i=0$ (and $\xi_i=0$),
and so, $\lambda_i'=\beta_i=\mu_i'$. Finally, if $a_i<x_i'<b_i$, then $\lambda_i'=0$
and $\mu_i'=1$, and so, since $\beta_i\in\{0,1\}$, then $\lambda_i'\le\beta_i\le\mu_i'$.

The uniqueness of $\beta(\alpha)$, for each $1-\xi\le\alpha\le1$, is a consequence
of the following: if $\lambda'\le\beta\le\mu'$ and $\beta\ne\beta(\alpha)$, then
there exists $i\in\{1,\dots,n\}$ s.t.\ $\beta_i=1-\beta_i(\alpha)$. Arguing as in
(a) and (b) above, we find that the equality \eq{e:E.9} cannot hold for this $\beta$.

Now, inequality \eq{e:E.4} readily follows from Lemma~\ref{l:Efour}, \eq{e:E.8}
and \eq{e:E.2}.
\qed\end{prEtwo}

\begin{rem} \label{r:abc}
(a) If $x'\in I_x^y$ in Lemma~\ref{l:Etwo}, then it is easily seen that $\xi=\mu-\lambda$,
and so, equality \eq{e:E.2} assumes the form:
  $$I_x^y=\bigcup_{1-(\mu-\lambda)\le\alpha\le1}I_{x+\alpha(\mu-\lambda)(x'-x)}^{x'+\alpha(y-x')}.$$
Although this equality looks different from \eq{e:E.6}, the two equalities are the same:
this is verified as in (i)--(iii) of the proof of Lemma~\ref{l:Efour}.

\smallbreak
(b) If $x<x'<y$, then $\xi=1$, $\lambda=0$ and $\mu=1$, and so, \eq{e:E.2}, \eq{e:E.6}
and \eq{e:E.4,5} are identical.

\smallbreak
(c) Here we consider a certain particular case of \eq{e:E.9} and establish conditions
on $x'$, under which $x'$ does not induce a (further) partition of $I_x^y$.
In view of \eq{e:E.9}, we have:
  \begin{equation} \label{e:E.10}
I_{x+\alpha\xi(x'-x)}^{x'+\alpha(y-x')}=I_x^y\quad\mbox{if and only if $\xi=0$ and $\alpha=1$,}
  \end{equation}
which is also equivalent to
  \begin{equation} \label{e:E.11}
a+\beta(x'-a)\le x\quad\mbox{and}\quad y\le x'+\beta(b-x')\quad\mbox{with}
\quad\beta=\beta(1).
  \end{equation}
Clearly, if $\xi=0$ and $\alpha=1$, then the left-hand side equality in \eq{e:E.10}
holds. Conversely, if the left-hand side equality in \eq{e:E.10} holds for some
$\mbox{$1\!-\!\xi\!\le\!\alpha\!\le\!1$}$, then $x+\alpha\xi(x'-x)=x$ and $x'+\alpha(y-x')=y$,
and so, if we suppose that $\xi_i=1$ for some $i\in\{1,\dots,n\}$, then, by \eq{e:E.3},
$x_i<x_i'<y_i$, and so, $\alpha_i=0$ and $x_i'=y_i$, which is a contradiction.
Thus, $\xi=0$ and $\alpha=1$.

Now, if $\xi=0$ and $\alpha=1$, then, by \eq{e:E.9} and \eq{e:E.10},
  \begin{equation} \label{e:E.12}
I_x^y=I_x^y\cap I_{a+\beta(x'-a)}^{x'+\beta(b-x')}\quad\mbox{with}\quad\beta=\beta(1),
  \end{equation}
which implies \eq{e:E.11}. Conversely, \eq{e:E.11} implies \eq{e:E.12}, and so, the
left-hand side equality in \eq{e:E.10} holds, i.e., $\xi=0$ and $\alpha=1$.

This observation also shows that a point $x'\in I_a^b$ induces a `true' partition of
$I_x^y$ provided that, for all $\beta$ with $\lambda'\le\beta\le\mu'$, we have:
  $$a+\beta(x'-a)\not\le x\quad\mbox{or}\quad y\not\le x'+\beta(b-x'),$$
which is also equivalent to $\xi\ne0$.
\end{rem}

\begin{prEone}
Let $\mathcal{P}=\{x[\sigma]\}_{\sigma=0}^\kappa$ for some $\kappa\in\mathbb{N}^n$
and $x'\in\mathcal{P}'$. Given $1\le\sigma\le\kappa$, we set $x_\sigma=x[\sigma-1]$,
$y_\sigma=x[\sigma]$ and $\xi_\sigma(x')=\xi(x_\sigma,x',y_\sigma)$, where $\xi$
is defined in \eq{e:E.3}, and note that $x_\sigma<y_\sigma$. The point $x'$ induces
a partition of $I_{x_\sigma}^{y_\sigma}=I_{x[\sigma-1]}^{x[\sigma]}$ of the form
\eq{e:E.2} with $x=x_\sigma$ and $y=y_\sigma$, and so, by virtue of \eq{e:E.1},
we get the following partition of $I_a^b$, induced by $x'$:
  \begin{equation} \label{e:E.14}
I_a^b=\bigcup_{1\le\sigma\le\kappa}\,\,\bigcup_{1-\xi_\sigma(x')\le\alpha\le1}
I_{x_\sigma+\alpha\xi_\sigma(x')(x'-x_\sigma)}^{x'+\alpha(y_\sigma-x')}.
  \end{equation}
We denote by $\mathcal{P}^1$ the net partition of $I_a^b$ corresponding to \eq{e:E.14}.
Moreover, by \eq{e:E.4}, for each $1\le\sigma\le\kappa$ we have the inequality:
  \begin{equation} \label{e:E.15}
\mbox{md}_n(f,I_{x_\sigma}^{y_\sigma})\le\sum_{1-\xi_\sigma(x')\le\alpha\le1}
\mbox{md}_n\bigl(f,I_{x_\sigma+\alpha\xi_\sigma(x')(x'-x_\sigma)}^{x'+\alpha(y_\sigma-x')}\bigr).
  \end{equation}

With no loss of generality we may assume that $x'\notin\mathcal{P}$:
if $x'\in\mathcal{P}$, i.e., $x'=x[\sigma']$ for some $1\le\sigma'\le\kappa$,
then $x'$ does not affect the partition $\mathcal{P}$ of $I_a^b$ in the sense
that $\mathcal{P}^1=\mathcal{P}$, and so, $\mbox{v}_n(f;\mathcal{P}^1)=\mbox{v}_n(f;\mathcal{P})$.
In order to see this, we note that \eq{e:E.3} implies
$\xi_\sigma(x')=\xi(x[\sigma-1],x[\sigma'],x[\sigma])=0$, and so, by Remark~\ref{r:abc}(c),
conditions \eq{e:E.10} and \eq{e:E.11} hold with $\beta=\beta(1)=(\beta_1,\dots,\beta_n)$ s.t.\
  $$\beta_i=\left\{
    \begin{array}{ccl}
    1 &\mbox{if}& x_i(\sigma_i')<x_i(\sigma_i),\\[2pt]
    0 &\mbox{if}& x_i(\sigma_i')\ge x_i(\sigma_i),
    \end{array}\right.
   =\left\{
    \begin{array}{ccl}
    1 &\mbox{if}& \sigma_i'<\sigma_i,\\[2pt]
    0 &\mbox{if}& \sigma_i'\ge \sigma_i.
    \end{array}\right.$$

Summing over $1\le\sigma\le\kappa$ in \eq{e:E.15} and taking into account \eq{e:E.1}
and \eq{e:E.14}, we obtain the inequality
  $$\mbox{v}_n(f;\mathcal{P})\le\mbox{v}_n(f;\mathcal{P}^1).$$

Replacing $\mathcal{P}$ by $\mathcal{P}^1$ in the arguments above, taking
$x'\in\mathcal{P}'\setminus\mathcal{P}^1$ and denoting by $\mathcal{P}^2$
the partition of $I_a^b$ induced from $\PC^1$ by $x'$, we get
$\mbox{v}_n(f;\PC^1)\le\mbox{v}_n(f;\PC^2)$. Since $\PC'\setminus\PC$ is a finite
set, we exhaust it by points $x'$ in a finite number of steps, arrive at the
partition $\PC'$ of $I_a^b$ and prove the desired inequality
$\mbox{v}_n(f;\PC)\le\mbox{v}_n(f;\PC')$.
\qed\end{prEone}

\begin{prLtroix}
1. First, we establish \eq{e:addi} for $\alpha=1=1_n$, i.e.,
  \begin{equation} \label{e:E.17}
V_n(f,I_x^y)=\sum_{1\le\sigma\le\kappa}V_n(f,I_{x[\sigma-1]}^{x[\sigma]}).
  \end{equation}
Modulo the notation, there is no loss of generality if we assume that $x=a$ and $y=b$,
so that $\{x[\sigma]\}_{\sigma=0}^\kappa$ is a net partition of $I_a^b$.

Let $\PC$ be an arbitrary net partition of $I_a^b$. Denote by $\PC'$ the net
partition of $I_a^b$ induced from $\PC$ by points $\{x[\sigma]\}_{\sigma=0}^\kappa$,
so that $\PC'$ is a refinement of $\PC$. Given $1\le\sigma\le\kappa$, set
$\PC_\sigma=\PC'\cap I_{x[\sigma-1]}^{x[\sigma]}$ and note that $\PC_\sigma$
is a net partition of $I_{x[\sigma-1]}^{x[\sigma]}$, and
$\PC'=\bigcup_{1\le\sigma\le\kappa}\PC_\sigma$. Then by virtue of Lemma~\ref{l:Eone},
we have:
  $$\mbox{v}_n(f;\PC)\le\mbox{v}_n(f;\PC')=\sum_{1\le\sigma\le\kappa}\mbox{v}_n(f;\PC_\sigma)
    \le\sum_{1\le\sigma\le\kappa}V_n(f,I_{x[\sigma-1]}^{x[\sigma]}).$$
Since $\PC$ is arbitrary, the left-hand side in \eq{e:E.17} is not greater than
the right-hand side.

Let us prove the reverse inequality. If $V_n(f,I_{x[\sigma-1]}^{x[\sigma]})$ is infinite
for some $1\le\sigma\le\kappa$, then since $I_{x[\sigma-1]}^{x[\sigma]}\subset I_a^b=I_x^y$,
the value $V_n(f,I_x^y)$ is infinite as well. Thus, we suppose that the right-hand side
of \eq{e:E.17} is finite. Let $\varepsilon>0$ be arbitrary. Given
$1\le\sigma\le\kappa$, by the definition of $V_n(f,I_{x[\sigma-1]}^{x[\sigma]})$,
there exists a net partition of $I_{x[\sigma-1]}^{x[\sigma]}$, denoted by
$\PC_\sigma(\varepsilon)$, s.t.\
  $$\mbox{v}_n(f;\PC_\sigma(\varepsilon))\ge V_n(f,I_{x[\sigma-1]}^{x[\sigma]})
    -(\varepsilon/c),$$
where $c=|\{\sigma:1\le\sigma\le\kappa\}|$. We denote by $\PC(\varepsilon)$ the net
partition of $I_a^b$ induced from $\{x[\sigma]\}_{\sigma=0}^\kappa$ by points
from $\bigcup_{1\le\sigma\le\kappa}\PC_\sigma(\varepsilon)$. Given $1\le\sigma\le\kappa$,
we set $\PC_\sigma'(\varepsilon)=\PC(\varepsilon)\cap I_{x[\sigma-1]}^{x[\sigma]}$ and
note that $\PC_\sigma'(\varepsilon)$ is a refinement of $\PC_\sigma(\varepsilon)$,
and $\PC(\varepsilon)=\bigcup_{1\le\sigma\le\kappa}\PC_\sigma'(\varepsilon)$.
By virtue of Lemma~\ref{l:Eone}, we find
  \begin{eqnarray}
V_n(f,I_a^b)&\ge&\mbox{v}_n(f;\PC(\varepsilon))=\sum_{1\le\sigma\le\kappa}
  \mbox{v}_n(f;\PC_\sigma'(\varepsilon))\ge\sum_{1\le\sigma\le\kappa}
  \mbox{v}_n(f;\PC_\sigma(\varepsilon))\nonumber\\
&\ge&\sum_{1\le\sigma\le\kappa}V_n(f,I_{x[\sigma-1]}^{x[\sigma]})
  -\varepsilon\biggl(\sum_{1\le\sigma\le\kappa}1\biggr)\Bigl/c,\nonumber
  \end{eqnarray}
where the factor by $\varepsilon$ is, actually, equal to $1$. The desired inequality
follows if we take into account the arbitrariness of $\varepsilon>0$.

\smallbreak
2. Now, suppose that $0\ne\alpha\le1$ and $\alpha\ne1$. Note that
$x\lfloor\alpha,y\lfloor\alpha\in I_{a\lfloor\alpha}^{b\lfloor\alpha}$ and
$x\lfloor\alpha<y\lfloor\alpha$, and that $\{x[\sigma]\lfloor\alpha\}_{\sigma\lfloor\alpha=0}%
^{\kappa\lfloor\alpha}$ is a net partition of $I_{a\lfloor\alpha}^{b\lfloor\alpha}$.
So, replacing $1=1_n$ by $1\lfloor\alpha$ (so that $|1\lfloor\alpha|=|\alpha|$)
and $f$---by $f_\alpha^z$ in \eq{e:E.17}, we get:
  \begin{eqnarray}
V_{|\alpha|}(f_\alpha^z,I_x^y\lfloor\alpha)&=&V_{|1\lfloor\alpha|}(f_\alpha^z,
  I_{x\lfloor\alpha}^{y\lfloor\alpha})\nonumber\\[2pt]
&=&\sum_{1\lfloor\alpha\le\sigma\lfloor\alpha\le\kappa\lfloor\alpha}
  V_{|1\lfloor\alpha|}\bigl(f_\alpha^z,I_{x[\sigma-1]\lfloor\alpha}^{x[\sigma]\lfloor\alpha}\bigr),
  \nonumber
  \end{eqnarray}
which is equal to the right-hand side of \eq{e:addi}.

This completes the proof of Lemma~\ref{l:troix}.
\qed\end{prLtroix}

\section{Proof of Theorem \ref{t:E}} \label{s:prE}

\begin{prtE}
1. First, we show that if $x,y\in I_a^b$, $x<y$, and $0\ne\alpha\le1$, then
  \begin{equation} \label{e:E*.1}
\mbox{md}_{|\alpha|}(f_\alpha^a,I_x^y\lfloor\alpha)=
\lim_{j\to\infty}\,\mbox{md}_{|\alpha|}((f_j)_\alpha^a,I_x^y\lfloor\alpha).
  \end{equation}
By virtue of \eq{e:2.4*}, we have:
  $$\mbox{md}_{|\alpha|}(f_\alpha^a,I_x^y\lfloor\alpha)=d\biggl(%
    \sum_{\even\theta\le\alpha}f(\,\underbrace{a+\alpha(x-a)+\theta(y-x)}_{({\displaystyle\cdots})}\,),%
    \sum_{\odd\theta\le\alpha}f(\cdots)\biggr),$$
and a similar equality holds for $f_j$ in place of $f$. Applying the inequalities
$|d(u,v)-d(u',v')|\le d(u,u')+d(v,v')$, $u,v,u',v'\in M$, and \eq{e:ms.2} and
taking into account the pointwise convergence of $f_j$ to $f$, we find
  \begin{eqnarray}
&&\bigl|\mbox{md}_{|\alpha|}((f_j)_\alpha^a,I_x^y\lfloor\alpha)-
  \mbox{md}_{|\alpha|}(f_\alpha^a,I_x^y\lfloor\alpha)\bigr|\nonumber\\[2pt]
&&\quad\,\,\le d\biggl(\sum_{\even\theta\le\alpha}f_j(\cdots),\sum_{\even\theta\le\alpha}f(\cdots)\biggr)
  +d\biggl(\sum_{\odd\theta\le\alpha}f_j(\cdots),\sum_{\odd\theta\le\alpha}f(\cdots)\biggr)\nonumber\\
&&\quad\,\,\le\sum_{\even\theta\le\alpha}d(f_j(\cdots),f(\cdots))+
  \sum_{\odd\theta\le\alpha}d(f_j(\cdots),f(\cdots))\nonumber\\
&&\quad\,\,=\sum_{0\le\theta\le\alpha}d(f_j(\cdots),f(\cdots))\to0\quad\mbox{as}\quad j\to\infty.\nonumber
  \end{eqnarray}

\smallbreak\label{pg:li}
2. In the rest of this proof we need only the inequality
  \begin{equation} \label{e:almin}
\mbox{md}_{|\alpha|}(f_\alpha^a,I_x^y\lfloor\alpha)\le
\liminf_{j\to\infty}\,\mbox{md}_{|\alpha|}((f_j)_\alpha^a,I_x^y\lfloor\alpha),
  \end{equation}
which readily follows from \eq{e:E*.1} and is applied one more time in the proof
of Theorem~\ref{t:weak} (Step~5). If $\mathcal{P}=\{x[\sigma]\}_{\sigma=0}^\kappa$ is a net
partition of $I_a^b$, then $\mathcal{P}\lfloor\alpha=\{x[\sigma]\lfloor\alpha\}%
_{\sigma\lfloor\alpha}^{\kappa\lfloor\alpha}$ is net partition of $I_a^b\lfloor\alpha$,
and so, given $1\le\sigma\le\kappa$, setting $x=x[\sigma-1]$ and $y=x[\sigma]$
in \eq{e:almin}, we find
  \begin{eqnarray}
\sum_{1\lfloor\alpha\le\sigma\lfloor\alpha\le\kappa\lfloor\alpha}
  \mbox{md}_{|\alpha|}(f_\alpha^a,I_{x[\sigma-1]}^{x[\sigma]}\lfloor\alpha)&\le&
  \sum_{1\lfloor\alpha\le\sigma\lfloor\alpha\le\kappa\lfloor\alpha}
  \liminf_{j\to\infty}\,\mbox{md}_{|\alpha|}((f_j)_\alpha^a,
  I_{x[\sigma-1]}^{x[\sigma]}\lfloor\alpha)\nonumber\\
&\le&\liminf_{j\to\infty}\sum_{1\lfloor\alpha\le\sigma\lfloor\alpha\le\kappa\lfloor\alpha}
  \mbox{md}_{|\alpha|}((f_j)_\alpha^a,I_{x[\sigma-1]}^{x[\sigma]}\lfloor\alpha)\nonumber\\
&\le&\liminf_{j\to\infty}\,V_{|\alpha|}((f_j)_\alpha^a,
  I_{x[\sigma-1]}^{x[\sigma]}\lfloor\alpha).\nonumber
  \end{eqnarray}
By the arbitrariness of $\mathcal{P}$, we infer that
  $$V_{|\alpha|}(f_\alpha^a,I_{x[\sigma-1]}^{x[\sigma]}\lfloor\alpha)\le
    \liminf_{j\to\infty}\,V_{|\alpha|}((f_j)_\alpha^a,I_{x[\sigma-1]}^{x[\sigma]}\lfloor\alpha).$$
We conclude that
  \begin{eqnarray}
\mbox{TV}(f,I_a^b)&=&\sum_{0\ne\alpha\le1}V_{|\alpha|}(f_\alpha^a,I_a^b\lfloor\alpha)
  \le\sum_{0\ne\alpha\le1}\liminf_{j\to\infty}\,V_{|\alpha|}((f_j)_\alpha^a,I_a^b\lfloor\alpha)\nonumber\\
&\le&\liminf_{j\to\infty}\,\sum_{0\ne\alpha\le1}V_{|\alpha|}((f_j)_\alpha^a,I_a^b\lfloor\alpha)
  =\liminf_{j\to\infty}\,\mbox{TV}(f_j,I_a^b).\nonumber\qquad\qquad\square
  \end{eqnarray}
\end{prtE}

\end{document}